\documentclass
[
    a4paper,
    DIV=11,
    abstracton,
]
{scrartcl}

\usepackage
{
    graphicx,
    amssymb,
    amsmath,
    amsthm,
    dsfont, 
    xcolor,
    mathtools,
    kbordermatrix,
    authblk,
    enumitem,
    tikz,
    todonotes,
    algorithm,
    algpseudocode,
    booktabs,
    stmaryrd, 
}

\usepackage[pdffitwindow=false,
            plainpages=false,
            pdfpagelabels=true,
            pdfpagemode=UseOutlines,
            pdfpagelayout=SinglePage,
            bookmarks=false,
            colorlinks=true,
            hyperfootnotes=false,
            linkcolor=blue,
            urlcolor=blue!30!black,
            citecolor=green!50!black]{hyperref}

\usepackage[bf,normal]{caption}
\usepackage{subcaption}

\graphicspath{{figures/}}

\newcommand{\R}{\mathbb{R}}                                     
\newcommand{\pd}[2]{\frac{\partial#1}{\partial#2}}              
\newcommand{\pdtwo}[3]{\frac{\partial^2 #1}{\partial #2 \partial #3}}			
\newcommand{\innerprod}[2]{\left\langle #1,\, #2 \right\rangle} 
\newcommand{\ts}{\hspace*{0.1em}}                               
\newcommand{\PSI}{\boldsymbol{\Psi}}
\DeclareMathOperator{\spn}{span}

\newcommand{\T}{\mathbf{T}}
\newcommand{\p}{\psi}
\newcommand{\e}{\mathbf{e}}
\newcommand{\f}{\mathbf{f}}
\newcommand{\g}{\mathbf{g}}
\newcommand{\h}{\mathbf{h}}
\newcommand{\A}{\mathbf{A}}
\newcommand{\B}{\mathbf{B}}

\newcommand{\super}[1]{^{(#1)}}  
\newcommand{\core}[1]{\left\llbracket #1 \right\rrbracket}  

\newtheorem{theorem}{Theorem}[section]
\newtheorem{corollary}[theorem]{Corollary}
\newtheorem{lemma}[theorem]{Lemma}
\newtheorem{proposition}[theorem]{Proposition}
\newtheorem{definition}[theorem]{Definition}
\theoremstyle{definition}

\newtheorem{remark}[theorem]{Remark}

\mathtoolsset{centercolon} 

\allowdisplaybreaks

\makeatletter
\def\blfootnote{\gdef\@thefnmark{}\@footnotetext}
\makeatother

\title{tgEDMD: Approximation of the Kolmogorov Operator in Tensor Train Format}
\author[1]{Marvin L\"ucke*}
\author[2, 3]{Feliks N\"uske*}
\affil[1]{Zuse Institute Berlin, Germany}
\affil[2]{Institute of Mathematics, Paderborn University, Germany}
\affil[3]{Max Planck Institute for Dynamics of Complex Technical Systems, Magdeburg, Germany}

\date{\today}

\begin{document}
\maketitle
\blfootnote{* Authors contributed equally to this work}
\blfootnote{Correspondence: nueske@mpi-magdeburg.mpg.de}

\begin{abstract}
Extracting information about dynamical systems from models learned off simulation data has become an increasingly important research topic in the natural and engineering sciences. Modeling the Koopman operator semigroup has played a central role in this context. As the approximation quality of any such model critically depends on the basis set, recent work has focused on deriving data-efficient representations of the Koopman operator in low-rank tensor formats, enabling the use of powerful  model classes while avoiding over-fitting. On the other hand, detailed information about the system at hand can be extracted from models for the infinitesimal generator, also called Kolmogorov backward operator for stochastic differential equations. In this work, we present a data-driven method to efficiently approximate the generator using the tensor train (TT) format. The centerpiece of the method is a TT representation of the tensor of generator evaluations at all data sites. We analyze consistency and complexity of the method, present extensions to practically relevant settings, and demonstrate its applicability to benchmark numerical examples.
\end{abstract}

\section{Introduction}
Learning models for complex dynamical systems based on simulation or measurement data has become a vibrant and successful research field, with applications in fluid dynamics, engineering, biophysics, economics, and many others. A significant body of work revolves around the statistical description of a dynamical system using the \emph{Koopman operator}, see~\cite{KOOPMAN1931,Dellnitz1999,Schuette1999,Mezic2005,Budisic2012,KLUS2016numerical,Klus2018a,MAUROY2020}.
Major algorithmic frameworks originating from this line of research include, for example, \emph{(extended) dynamic mode decomposition}, or (E)DMD~\cite{SCHMID2010,WILLIAMS2015}, \emph{Markov state modeling}~\cite{Dellnitz1999,Schuette1999,Prinz2011c,BowmanPandeNoe_MSMBook}, the \emph{variational approach to conformational dynamics} (VAC)~\cite{Noe2013,Nueske2014} and its generalization \emph{variational approach to Markov processes} (VAMP)~\cite{WU2020}, and the SINDy approach to system identification~\cite{BRUNTO2016}. At their core, most of these techniques aim at solving a linear regression problem, the solution of which consistently approximates the Koopman operator on a pre-selected model class in the infinite data limit. In order to reduce the influence of statistical noise when using finite data, the problem dimension is often reduced via a projection onto the leading components of a singular value decomposition (SVD) of the data matrix. This de-noising and whitening technique is also known as \emph{AMUSE} algorithm~\cite{TONG1990,Klus2018a}. Non-linear machine learning techniques have also been incorporated into the Koopman framework, see~\cite{MARDT2018,Lusch2018,Otto2019} for deep learning approaches, and~\cite{Williams2015kernel,KLUS2020eigen,KLUS2020kernel} for kernel-based methods.

\noindent This work combines two recent developments along these lines: the first is the use of a specific low-rank tensor representation, the \emph{tensor train format} (TT format), in order to balance sufficient representational power of the model class and manageable computational complexity. The TT format has been known in various scientific fields for decades, especially in quantum chemistry~\cite{AFFLECK1987,OSTLUND1995,SZALAY2015}. It was introduced in the mathematical literature in~\cite{Oseledets2009,Oseledets2011}, and was applied specifically to the Fokker-Plack equation of an SDE in~\cite{Chertkov2021}. Within the Koopman framework, it was first introduced in~\cite{NUESKE2016,KLUS2016}. Building on results from~\cite{KLUS2018}, a multi-linear version of the AMUSE algorithm to be used with the TT format was presented by one of the authors in~\cite{Nueske2021}. In particular, it was shown that, due to the TT structure of all quantities involved, one can efficiently compute a representation of the Koopman operator on a reduced subspace defined by a multi-linear SVD. Consistency of this representation for fixed TT ranks was also proven.

\noindent The second development we will build on is the extension of EDMD and related approaches to approximate the infinitesimal generator of the Koopman operator semigroup. The EDMD method can be adapted in a fairly direct manner for the purpose of generator approximation, as shown in~\cite{KLUS2020}. The resulting method, called \emph{gEDMD}, can further be utilized to identify system parameters, to define coarse grained models, and for system control. For alternative approaches to generator approximation, please see~\cite{FrJuKo13,Kaiser2018,Gia19,Mauroy2019,Kaiser2021}.

\noindent The core contribution of this paper is to combine these developments for systems driven by stochastic differential equations (SDE), where the Koopman generator is also called Kolmogorov backward operator. The result is a data-driven method to approximate the Kolmogorov operator on a multi-linear subspace of a tensor product function space. In more detail, we
\begin{enumerate}
\item derive a TT representation formula for tensors which are given as a certain structured sum of rank-one tensors, extending a result from~\cite{Kazeev2013}. This representation is then used to deduce the TT structure of the tensor of generator evaluations at all data sites, which is the missing piece when translating existing methods to the approximation of the generator. We derive separate representations for reversible and general SDE systems.
\item present the resulting data-driven algorithm, called tgEDMD, which combines the methods suggested in Refs.~\cite{KLUS2020} and ~\cite{Nueske2021}. Further, we analyze its computational complexity and consistency.
\item present extensions of the method to two practically relevant settings: first, we show how importance sampling ratios can be incorporated if the sampling distribution differs from the reference distribution for generator approximation. Second, we explain how basis sets defined on reduced variables can be used to approximate the generator of an effective dynamics in TT format, building on results from Refs.~\cite{Zhang2016,KLUS2020}.
\item illustrate the capabilities of the proposed methods by approximating the generator of a four-dimensional model system, and the effective generator for molecular dynamics of the deca-alanine peptide on a set of torsion angle coordinates.
\end{enumerate}

\noindent The rest of the paper is structured as follows: basic concepts regarding the Koopman framework, data-driven approximation, the AMUSE algorithm, and tensor trains are introduced in Section~\ref{sec:prerequisites}. The general TT representation formula and the TT structures of the generator data tensor are derived in Section~\ref{sec:tensor_reps}. The main algorithm, its analysis, and the extensions are presented in Section~\ref{sec:algorithms}. Numerical examples follow in Section~\ref{sec:numerical_results}.

\section{Prerequisites}
\label{sec:prerequisites}
\subsection{Kolmogorov Backward Operator}
\label{subsec:kolmogorov_op}
In this paper, we are concerned with the numerical approximation of the \emph{Kolmogorov backward operator}
\begin{equation}
\label{eq:def_generator}
\mathcal{L}\phi(x) = \frac{1}{2}a(x) : \nabla^2\phi(x) + b(x) \cdot \nabla\phi(x),
\end{equation}
where $b \in \mathbb{R}^d$ is a vector field and $a \in \mathbb{R}^{d \times d}$ is a field of symmetric positive semi-definite matrices (that is, $v^T a(x) v \geq 0$ for all $x$ and $v \in \mathbb{R}^d$). Moreover, $\nabla \phi(x) \in \mathbb{R}^d$ is the gradient of the scalar function $\phi$, $\nabla^2 \phi(x)$ is the Hessian matrix of its second-order derivatives, and the colon $:$ denotes the Frobenius inner product, $A:B = \mathrm{tr}(A^T B)$, between matrices. 

\noindent The differential operator $\mathcal{L}$ is closely linked to the statistical description of dynamical systems which are driven by a stochastic differential equation (SDE). For $t \geq 0$, let $X_t$ be a stochastic process on a state space $\mathbb{X} \subset \mathbb{R}^d$, subject to the governing equation
\begin{equation}
\label{eq:sde}
dX_t = b(X_t) \mathrm{d}t + \sigma(X_t) \mathrm{d}W_t.
\end{equation}
In this context, the field $b$ is called the \emph{drift}, while the matrix field $\sigma$ is called \emph{diffusion}, which must be such that $\sigma \sigma^T = a$. Furthermore, $W_t$ in~\eqref{eq:sde} is the Brownian motion in $d$-dimensional space. Equation~\eqref{eq:sde} can be understood in the sense that the infinitesimal change of the process $X_t$ is proportional to the drift field $b$, but also subject to the random fluctuations of the Brownian motion $W_t$, modulated by $\sigma$. See Ref.~\cite{Oksendal2013} for a rigorous introduction to the dynamics of type \eqref{eq:sde}.

\noindent The statistics of the process $X_t$ are then captured by a family of linear operators, called the \emph{Koopman semigroup} \cite{KOOPMAN1931,Mezic2005}. For $t \geq 0$ and a function $\phi$ on $\mathbb{X}$ (usually called an \emph{observable} of the system), the linear operator $\mathcal{K}^t$ acts as
\begin{equation*}
\mathcal{K}^t \phi(x) = \mathbb{E}^x[\phi(X_t)] = \mathbb{E}[\phi(X_t) \vert X_0 = x].
\end{equation*} 
The Kolmogorov operator $\mathcal{L}$ arises as the \emph{infinitesimal generator} of this semigroup, that is, for suitable functions $\phi$, we have
\begin{equation}
\label{eq:derivative_generator}
 \mathcal{L}\phi = \lim_{t \rightarrow 0} \frac{1}{t}(\mathcal{K}^t - \mathrm{Id})\phi.
\end{equation}
As described in the introduction, knowledge of the Koopman semigroup, and especially its generator, can be utilized for a thorough statistical analysis of the dynamical system $X_t$. Before discussing its numerical approximation, let us point out that the precise mathematical properties of the differential operator $\mathcal{L}$ and its domain of definition can be quite intricate. As these technical details are mostly irrelevant for the present study, we just briefly lay out a standard setting to serve as a guideline for the following, see~\cite{Oksendal2013,Bakry2014} for an in-depth discussion of this setting.

\begin{itemize}
\item[$\mathbf{(A1)}$] The coefficients $b$, $\sigma$ of the SDE~\eqref{eq:sde} are such that global existence and uniqueness of solutions with continuous paths are guaranteed.
\item[$\mathbf{(A2)}$] There is a unique, everywhere positive invariant density $\rho$ for~\eqref{eq:sde}, with associated probability measure $\mu$. Invariance of the measure $\mu$ can be characterized by the identity, for all $\phi$ and $t \geq 0$:
\begin{equation*}
\int_\mathbb{X} \mathcal{K}^t \phi(x) \, \mathrm{d}\mu(x) = \int_\mathbb{X} \phi(x) \, \mathrm{d}\mu(x).
\end{equation*}
\end{itemize}
It follows that the operators $\mathcal{K}^t$ form a strongly continuous contraction semigroup on the Hilbert space of square-integrable functions with respect to the invariant measure $\mu$:
\[ L^2_\mu(\mathbb{X}) = \{\phi:\, \mathbb{X}\mapsto \mathbb{R} \quad \vert \int_\mathbb{X} |\phi(x)|^2 \,\mathrm{d}\mu(x) < \infty \}. \]
The inner product for $\phi,\, \tilde{\phi} \in L^2_\mu(\mathbb{X})$ is given by $\innerprod{\phi}{\tilde{\phi}}_\mu = \int_\mathbb{X} \phi \tilde{\phi}\,\mathrm{d}\mu$. The domain $\mathcal{D}(\mathcal{L})$ of functions where the derivative \eqref{eq:derivative_generator} exists is a dense subspace of $L^2_\mu(\mathbb{X})$. If in addition the dynamics \eqref{eq:sde} are \emph{reversible} with respect to the invariant measure $\mu$, that is the detailed balance condition
\[ \mathbb{P}^\mu(X_0 \in A,\, X_t \in B) = \mathbb{P}^\mu(X_0 \in B,\, X_t \in A) \]
holds for all $A, B \subset \mathbb{X}$ and $t \geq 0$, then the generator $\mathcal{L}$ is self-adjoint on its domain. Moreover, the following important relation holds for all $\phi, \, \tilde{\phi} \in \mathcal{D}(\mathcal{L})$:

\begin{equation}
\label{eq:reversible_bilinear_form}
\innerprod{\mathcal{L \phi}}{\tilde{\phi}}_\mu = -\frac{1}{2}\int_\mathbb{X} \nabla \phi(x)^T \,  a(x) \, \nabla \tilde{\phi}(x) \,\mathrm{d}\mu(x).
\end{equation}

\subsection{Galerkin Projection and Data-driven Estimation}
\label{subsec:galerkin_projection}
Next, we consider the approximation of $\mathcal{L}$ on a finite-dimensional subspace $\mathbb{V}$, spanned by twice continuously differentiable functions contained in the domain, i.e. $\mathbb{V} \subset \mathcal{D}(\mathcal{L}) \cap C^2(\mathbb{X})$. Choosing $\mathbb{V}$ in this way ensures that $\mathcal{L}\phi(x)$ can be evaluated in the form~\eqref{eq:def_generator} for all $\phi \in \mathbb{V}$. If a basis set $\psi = [\psi_1, \ldots, \psi_n]^T$ for $\mathbb{V}$ is given, we will briefly write $\psi(x) = [\psi_1(x), \ldots, \psi_n(x)]^T \in \mathbb{R}^n$ for the evaluation of all basis functions at a point $x$, and similarly $\mathcal{L}\psi(x) \in \mathbb{R}^n$ for the application of $\mathcal{L}$ to all basis functions at once. Moreover, if $X = [x_1, \ldots, x_m] \in \mathbb{R}^{d \times m}$ is a collection of data points in $\mathbb{X}$, we write
\[\Psi(X) = [\psi(x_1),\, \ldots, \, \psi(x_m)] \in \mathbb{R}^{n \times m}, \quad \mathcal{L}\Psi(X) = [\mathcal{L}\psi(x_1),\, \ldots, \, \mathcal{L}\psi(x_m)] \in \mathbb{R}^{n\times m}, \]
and call these the \emph{data matrix} and \emph{generator data matrix}, respectively.

\noindent The orthogonal projection of $\mathcal{L}$ onto $\mathbb{V}$ can be represented by the matrix (w.r.t. the basis $\psi$):
\begin{align}
\nonumber L_\mathbb{V} &:= C(\psi)^{-1} A(\psi),  \\
\label{eq:galerkin_L}
C(\psi)_{uv} &= \innerprod{\psi_u}{\psi_v}_\mu = \mathbb{E}^\mu[\psi_u \ts \psi_v], & A(\psi)_{uv} &= \innerprod{\psi_u}{\mathcal{L}\psi_v}_\mu = \mathbb{E}^\mu[\psi_u \ts \mathcal{L}\psi_v].
\end{align}
As the integrals appearing in~\eqref{eq:galerkin_L} are expectation values with respect to the probability measure $\mu$, they can be approximated by Monte Carlo methods. From now on, we assume that $x_l, \, l\in \mathbb{N}$, is a sequence of random variables, such that for all integrable functions $\phi \in L^1_\mu(\mathbb{X})$ and almost all realizations of the $x_l$ we have: 
\begin{equation}
\label{eq:ergodic_limit}
\lim_{m \rightarrow \infty} \frac{1}{m}\sum_{l=1}^m \phi(x_l) = \int_\mathbb{X} \phi(x)\,\mathrm{d}\mu(x).
\end{equation}
Equation~\eqref{eq:ergodic_limit} is of course satisfied if all $x_l$ are i.i.d. with respect to $\mu$.  By the ergodic theorem~\cite{Chacon1962}, the same conclusion holds if $x_l$ is taken as the $l$-th step of any discretized trajectory of the dynamics~\eqref{eq:sde}, provided the initial condition is drawn from $\mu$. As a consequence, for any finite realization $x_1, \ldots, x_m$ of the random variables $x_l$, we obtain an empirical Galerkin matrix $\hat{L}_\mathbb{V} = \hat{C}(\psi)^{-1}\hat{A}(\psi)$, where

\begin{align}
\label{eq:empirical_estimator_C}
\hat{C}(\psi) &= \frac{1}{m}\sum_{l=1}^m \psi(x_l) \otimes \psi(x_l) = \frac{1}{m}\Psi(X)\Psi(X)^T, \\
\label{eq:empirical_estimator_A}
\hat{A}(\psi) &= \frac{1}{m}\sum_{l=1}^m \psi(x_l) \otimes \mathcal{L}\psi(x_l) = \frac{1}{m}\Psi(X)(\mathcal{L}\Psi(X))^T.
\end{align}
For the reversible case, the matrix $\hat{A}(\psi)$ can be calculated in a simpler form based on~\eqref{eq:reversible_bilinear_form}, requiring only first order derivatives of the basis set. Define $\mathrm{d}\psi(x) =  [\nabla \psi (x) \sigma(x)]\in \mathbb{R}^{n \times d}$, using the Jacobian $\nabla \psi \in \mathbb{R}^{n \times d}$. Then, introduce a modified generator data matrix $\mathrm{d}\Psi(X) = \left[ \mathrm{d}\psi(x_1) \dots \mathrm{d}\psi(x_m)\right] \in \mathbb{R}^{n \times dm}$. It follows from~\eqref{eq:reversible_bilinear_form} that $A(\psi)$ can be approximated by
\begin{equation}
\label{eq:empirical_estimator_Arev}
\hat{A}(\psi) = -\frac{1}{2m} \mathrm{d}\Psi(X) [\mathrm{d}\Psi(X)]^T,
\end{equation}
noting that the summation extends over both the data size and the state space dimension. This conceptually simple way of approximating the Kolmogorov operator based on data using (\ref{eq:empirical_estimator_C} - \ref{eq:empirical_estimator_Arev}) has been labeled \emph{generator extended dynamic mode decomposition (gEDMD)} in \cite{KLUS2020}.

\subsection{AMUSE Algorithm}
\label{subsec:AMUSE}
The calculation of the Gramian matrix $\hat{C}(\psi)$ and its inverse may be avoided by means of the AMUSE Algorithm~\ref{alg:AMUSE} \cite{TONG1990,Klus2018a}, which also serves as a tool for dimension reduction.

\begin{algorithm}[htb]
  \caption{AMUSE}
  \label{alg:AMUSE}
  \setlength{\tabcolsep}{.5ex}
  \begin{tabular}{ll}
    \textbf{Input:} & data matrix $\Psi(X)$, generator data matrix $\mathcal{L}\Psi(X)$ or $\mathrm{d}\Psi(X)$. \\
    \textbf{Output:} & Reduced matrix representation $\hat{M}_r \in \mathbb{R}^{r \times r}$ or $\hat{M}_r^{\mathrm{rev}} \in \mathbb{R}^{r \times r}$ of $\mathcal{L}$, where $r \leq n$.
  \end{tabular}
  \hrule\vspace{0.2cm}
  \begin{algorithmic}[1]
    \State Compute rank-$r$ SVD of $\Psi(X)$, i.e., $\Psi(X) \approx \hat{U}_r \, \hat{\Sigma}_r \, \hat{V}_r^\top$. \label{algline:AMUSE_SVD}
    \State \textbf{(Non-reversible Case)}: Compute reduced matrix\label{algline: AMUSE_reduced_matrix}
    \Statex \quad\quad\quad $\hat{M}_r = \hat{\Sigma}_r^{-1} \, \hat{U}_r^\top \, \Psi(X) \, \mathcal{L}\Psi(X)^T \, \hat{U}_r \, \hat{\Sigma}_r^{-1} = \hat{V}_r^\top \, \mathcal{L}\Psi(X)^T \, \hat{U}_r \, \hat{\Sigma}_r^{-1}$.
    \State \textbf{(Reversible Case)}: Compute reduced matrix\label{algline: AMUSE_reduced_matrix_rev}
    \Statex \quad\quad\quad $\hat{M}_r^{\mathrm{rev}} = -\frac{1}{2} \hat{\Sigma}_r^{-1} \, \hat{U}_r^\top \, \mathrm{d}\Psi(X) \, \mathrm{d}\Psi(X)^T \, \hat{U}_r \, \hat{\Sigma}_r^{-1}$.
  \end{algorithmic}
\end{algorithm}

\noindent The rationale behind AMUSE is that, using the truncated SVD $\Psi(X) \approx \hat{U}_r \, \hat{\Sigma}_r \, \hat{V}_r^\top$ computed in line~\ref{algline:AMUSE_SVD}, the transformed basis $\hat{\eta}_r(x)^T := \psi(x)^T \hat{U}_r \hat{\Sigma}_r^{-1}$  is empirically orthonormal, hence computation of the Gramian matrix is no longer necessary. In order to compute a representation of the generator on the linear span of $\hat{\eta}_r$, the same transformation just needs to be applied to $\hat{A}(\psi)$, which is what happens in lines \ref{algline: AMUSE_reduced_matrix} and \ref{algline: AMUSE_reduced_matrix_rev} of Algorithm \ref{alg:AMUSE}. The reduced matrices $\hat{M}_r$ and $\hat{M}_r^{\mathrm{rev}}$ can also be computed by the following expressions, which do not require assembly of the full generator data matrices:

\begin{align}
\label{eq:reduced_matrix_pw_nonrev}
\hat{M}_r &= \hat{V}_r^\top \, \left[\sum_{l=1}^m \mathcal{L}\psi(x_l) \otimes e_l\right]^T \, \hat{U}_r \, \hat{\Sigma}_r^{-1} = \sum_{l=1}^m (\hat{V}_r)_{l, :} \otimes \mathcal{L}\psi(x_l)^T \hat{U}_r \, \hat{\Sigma}_r^{-1}, \\
\label{eq:reduced_matrix_pw_rev}
\hat{M}_r^{\mathrm{rev}} &= -\frac{1}{2}\sum_{l=1}^m \hat{\Sigma}_r^{-1} \, \hat{U}_r^\top \, \nabla \psi(x_l)\, a(x_l) \, \nabla \psi(x_l)^T \, \hat{U}_r \, \hat{\Sigma}_r^{-1}.
\end{align}

\noindent If one is interested in consistency of the reduced matrix for the infinite data limit, the spectral decomposition of the analytical Gramian needs to be considered. Denote this spectral decomposition by $C(\psi) = U\Sigma^2 U^T$, and for any $r \leq n$ we write $C(\psi) = U_r \Sigma_r^2 U_r^T + E_r$ for its truncation using $r$ terms, with error term $E_r$. Now, if the reduced dimension $r$ in Algorithm~\ref{alg:AMUSE} is fixed, it can be shown that, unless the eigenvalue $\sigma_r^2$ of $C(\psi)$ is degenerate, the linear span of $\hat{\eta}_r$ converges almost surely to that of $\eta_r := \psi(x)^T U_r \Sigma_r^{-1}$, along with the projected generator on this space (see \cite{Nueske2021}[Proposition 3] for a proof). Alternatively, one can start with a fixed truncation parameter $\epsilon > 0$. Assume there is an $r(\epsilon)$ such that $\sigma_{r(\epsilon)} > \epsilon > \sigma_{r(\epsilon) + 1}$ in the spectral decomposition of $C(\psi)$. Then, if the reduced rank in Algorithm~\ref{alg:AMUSE} is chosen such that $\hat{\sigma}_r \geq \sqrt{m} \epsilon > \hat{\sigma}_{r + 1}$, it follows by the same arguments that $r \rightarrow r(\epsilon)$ almost surely as $m\rightarrow \infty$, recalling that $\hat{C}(\psi) = \frac{1}{m}\Psi(X)\Psi(X)^T$. Hence, $\hat{M}_r, \, \hat{M}_r^{\mathrm{rev}}$ provide consistent estimates of the projected generator on the linear span of $\eta_{r(\epsilon)}$. This data-dependent truncation strategy will also be used in the multi-linear case, see the discussion in Section~\ref{subsec:efficient_calculation_M} below.

\subsection{Tensor Train Format}
The choice of the finite-dimensional subspace $\mathbb{V}$ in Section~\ref{subsec:galerkin_projection} is of critical importance to the quality of the gEDMD approximation. A widely used approach to arrive at a large space with rich approximation properties is to consider all products of functions contained in a selection of elementary subspaces of moderate dimension. Many of the quantities introduced above will be endowed with a \emph{tensor structure} in this case. For our purposes, it is sufficient to think of tensors as multi-dimensional arrays, that is, array entries are labeled by $p$ different indices:
\[
\mathbf{T} = \mathbf{T}_{u_1, \ldots, u_p},\quad 1 \leq u_k \leq n_k, \quad 1\leq k \leq p.
\]
The number of indices $p$ is called \emph{order} of the tensor, the individual indices are often called \emph{modes} with \emph{mode sizes} $n_k$. A tensor can be \emph{unfolded} into a matrix by grouping any selection $\mathrm{I} \subset \{1, \ldots, p\}$ of tensor indices along rows and all remaining tensor indices along columns. We denote this matrix by $\mathbf{T}\vert_{\mathrm{I}}$. Specifically, if $\mathrm{I}$ consists of the first $k$ indices, we just write $\mathbf{T}\vert_k $ and call the resulting matrix the \emph{mode-$k$ unfolding}.

\noindent The \emph{tensor product} of $p$ vectors $v^{(k)} \in \mathbb{R}^{n_k}$ is defined by
\[\mathbf{T} = v^{(1)} \otimes \ldots \otimes v^{(p)}, \quad \mathbf{T}_{u_1, \ldots, u_p} = \prod_{k=1}^p v^{(k)}_{u_k}, \]
generalizing the standard outer product between two vectors. A tensor of this type is a \emph{rank-one} tensor. Every tensor can be written as a linear combination of finitely many rank-one tensors, but their number may be exponentially large. However, many tensors appearing in applications can be (approximately) represented by a much smaller number of parameters. Many different parametrizations, typically called \emph{tensor formats}, have been put forward. Here, we focus exclusively on the tensor train format (TT format), introduced in~\cite{Oseledets2009,Oseledets2011}, where a tensor is represented as a contraction of multiple lower-order tensors:

\begin{definition}
  A tensor $\mathbf{T}$ is said to be in the \emph{TT format} if
  \begin{equation}
  \label{eq:def_TT}
    \mathbf{T} = \sum_{s_0=1}^{r_0} \cdots  \sum_{s_p=1}^{r_p} \bigotimes_{k=1}^{p} \mathbf{T}^{(k)}_{s_{k-1},:,s_k} = \sum_{s_0=1}^{r_0} \cdots  \sum_{s_p=1}^{r_p}  \mathbf{T}^{(1)}_{s_0,:,s_1} \otimes \dots \otimes  \mathbf{T}^{(p)}_{s_{p-1},:,s_p}.
  \end{equation}
  The tensors $\mathbf{T}^{(k)} \in \R^{r_{k-1} \times n_k \times r_k}$ of order 3 are called \emph{TT cores} and the numbers $r_k$ are called \emph{TT ranks}. It holds that $r_0 = r_p =1$ and $r_k \geq 1$ for $k=1, \dots, p-1$. A TT core is \emph{left-orthonormal} if
\begin{equation*}
  \left( \mathbf{T}^{(k)}\vert_2 \right)^\top \cdot \left(\mathbf{T}^{(k)}\vert_2 \right) = \mathrm{Id} \in \R^{r_k \times r_k}.
\end{equation*}
\end{definition}

\begin{remark}
\label{rem:multiple_TT}
Sometimes, we will also allow the first or final rank to be greater than one, i.e. $r_0 >1$ or $r_p >1$. In this case, we understand $\mathbf{T}$ to be encoding multiple tensors, one for each pair of values of $s_0$ and $s_p$ in~\eqref{eq:def_TT}.
\end{remark}

\noindent The following formal representation for tensor trains has proven particularly useful. For a given tensor train $\mathbf{T}$ with cores $\mathbf{T}^{(k)} \in \R^{r_{k-1} \times n_k \times r_k}$, a single core is written as a two-dimensional matrix containing vectors as elements
\begin{equation}
  \left\llbracket \mathbf{T}^{(k)} \right\rrbracket = 
  \left\llbracket
  \begin{matrix}
    & \mathbf{T}^{(k)}_{1,:,1} & \cdots & \mathbf{T}^{(k)}_{1,:,r_k} & \\
    & & & & \\
    & \vdots & \ddots & \vdots & \\
    & & & & \\
    & \mathbf{T}^{(k)}_{r_{k-1},:,1} & \cdots & \mathbf{T}^{(k)}_{r_{k-1},:,r_k} &
  \end{matrix}\right\rrbracket.
  \label{eq:rank_core notation}
\end{equation}
We then represent $\mathbf{T}$ by a formal matrix product, called \emph{rank-core product} in~\cite{Kazeev2012}: 
\[ \mathbf{T} = \left\llbracket \mathbf{T}^{(1)}\right\rrbracket \otimes \dots \otimes \left\llbracket \mathbf{T}^{(p)}\right\rrbracket. \]
In this representation, we compute the tensor products of the matrix elements -- which are vectors instead of scalars -- and then sum over the columns and rows, just as for an ordinary matrix product. It is readily verified that the result is the same as~\eqref{eq:def_TT}.

\noindent Further below, we will also require a tensor product for matrices, which is known as \emph{Kronecker product}, and denoted by the same symbol. For matrices $A \in \mathbb{R}^{m \times n},\, B \in \mathbb{R}^{r \times q}$, its block-wise definition is
\begin{equation}
\label{eq:kronecker_product}
A \otimes B \in \mathbb{R}^{mr \times nq}: \quad A \otimes B = \begin{bmatrix}
A_{11} B & \ldots & A_{1n} B \\
\vdots & \ddots & \vdots \\
A_{m1} B & \ldots & A_{mn} B
\end{bmatrix}.
\end{equation}
The Kronecker product is equivalent to the standard tensor product applied to vectorized versions of the matrices $A, \, B$, but using a non-standard ordering of their indices, see \cite{HACKBUSCH2012} for details.

\subsection{Basis Decomposition and Global SVD}
\label{subsec:basis_global_svd}
We are now prepared to formulate the tensor-based approximation framework for the Kolmogorov operator. Suppose that for $1 \leq k \leq p$, $\mathbb{V}^k = \spn\{\psi_{k, 1}, \dots, \psi_{k, n_k}\} \subset \mathcal{D}(\mathcal{L}) \cap C^2(\mathbb{X})$ is an $n_k$-dimensional subspace. Assume further the tensor product 
\[ \mathbb{V} = \bigotimes_{k=1}^p \mathbb{V}^k = \spn\{\psi_{1, u_1} \cdot \ldots \cdot \psi_{p, u_p} \mid \, 1 \leq u_k \leq n_k, \, 1 \leq k \leq p \}\]
is also contained in $\mathcal{D}(\mathcal{L})$. We denote the basis formed by all the products above by $\PSI$ (assuming their linear independence). The evaluation of $\PSI$ at $x \in \mathbb{X}$ is a rank-one tensor of order $p$, while the analogue of the data matrix $\Psi(X)$ is a tensor of order $p+1$, which we call the \emph{data tensor}:
\[\PSI(X)_{u_1, \ldots, u_p, l} = \psi_{1, u_1}(x_l) \cdot \ldots \cdot \psi_{p, u_p}(x_l). \]
Tensor-based counterparts $\mathcal{L}\PSI(x),\, \mathcal{L}\PSI(X)$, and $\mathrm{d}\PSI(x), \,\mathrm{d}\PSI(X)$ of the generator data matrices can be defined the same way. Unless the $n_k$ and $p$ are very small, these data tensors cannot even be held in memory, let alone be used to compute the Galerkin matrices in (\ref{eq:empirical_estimator_C} - \ref{eq:empirical_estimator_Arev}). Therefore, we follow the lines of the AMUSE approach in Section~\ref{subsec:AMUSE}, but notice that we cannot compute an SVD of the unfolded data tensor $\PSI(X)\vert_p$ either, necessitating further modifications.

\noindent As shown in \cite{GELSS2019}, the data tensor is a sum of $m$ rank-one tensors, and hence possesses a TT representation with all ranks equal to $m$, and the following TT cores (using the rank-core notation introduced in \eqref{eq:rank_core notation}):
\begin{align}
\label{eq:tt_decomp_psi_x}
\core{\PSI(X)\super{1}} &= \core{\begin{matrix} \psi_1(x_1) & \dots & \psi_1(x_m) \end{matrix}}, &
\core{\PSI(X)\super{p+1}}  &= \core{\begin{matrix}
e_1 \\ \vdots \\ e_m
\end{matrix}}, \\
\nonumber \core{\PSI(X)\super{k}} &= \core{\begin{matrix}
    \psi_k (x_1) & & 0 \\
    & \ddots & \\
    0 & & \psi_k (x_m)
    \end{matrix}}, \, 2 \leq k \leq p. &
\end{align}
In the final core above, $e_l$ denote the canonical unit vectors in $\mathbb{R}^m$. Moreover, by a sequence of SVDs applied to the cores of~\eqref{eq:tt_decomp_psi_x}, a multi-linear approximation of the true SVD for $\PSI(X)\vert_p$ can be computed, called global SVD \cite{KLUS2018}, replacing the first step in Algorithm~\ref{alg:AMUSE}. The result is an SVD-like decomposition
\begin{equation*}
\PSI(X)\vert_p \approx \hat{\mathbf{U}}\vert_p \hat{\Sigma} \hat{V}^T.
\end{equation*}
Here, $\hat{\mathbf{U}} \in \mathbb{R}^{n_1 \times \ldots \times n_p \times r_p}$ is a tensor train with TT ranks $r_1, \ldots, r_p$, where $r_p > 1$ is allowed (see Remark~\ref{rem:multiple_TT}), and left-orthonormal cores. Moreover, $\hat{\Sigma} \in \mathbb{R}^{r_p \times r_p}$ is diagonal and non-negative, and $\hat{V} \in \mathbb{R}^{m \times r_p}$ is orthogonal. As shown in \cite{Nueske2021}, the columns of $\hat{\mathbf{U}}\vert_p \hat{\Sigma}^{-1}$ encode an empirically orthonormal basis of an $r_p$-dimensional subspace of $\mathbb{V}$:
\[ \hat{\eta}_{r_p}(x)^T = (\PSI(x)\vert_p)^T \hat{\mathbf{U}}\vert_p \hat{\Sigma}^{-1}. \]
Analogous to the standard AMUSE algorithm, we can then compute an empirical representation of the generator on the linear span of $\hat{\eta}_{r_p}$, by forming the product (see (\ref{eq:reduced_matrix_pw_nonrev} - \ref{eq:reduced_matrix_pw_rev})):
\begin{align}
\label{eq:amuse_approx_generator}
\hat{M}_{r_p} &= \sum_{l=1}^m \hat{V}_{l, :} \otimes (\mathcal{L}\PSI(x_l)\vert_p)^T \hat{\mathbf{U}}\vert_p \hat{\Sigma}^{-1}, \\
\label{eq:amuse_approx_generator_rev}
\hat{M}_{r_p}^{\mathrm{rev}} &= -\frac{1}{2} \sum_{l=1}^m \hat{\Sigma}^{-1} \, (\hat{\mathbf{U}}\vert_p)^T \, \nabla \PSI(x_l) \,a(x_l) \, \nabla \PSI(x_l)^T \, \hat{\mathbf{U}}\vert_p \, \hat{\Sigma}^{-1}.
\end{align}
The computationally demanding part of (\ref{eq:amuse_approx_generator}-\ref{eq:amuse_approx_generator_rev}) is the multiplication of the unfolded tensors $(\mathcal{L}\PSI(x_l)\vert_p)^T $ and $\hat{\mathbf{U}}\vert_p$ (or $(\nabla \PSI(x_l)\vert_p)^T $ and $\hat{\mathbf{U}}\vert_p$). The orthonormal factor $\hat{\mathbf{U}}$ is always given in tensor train format by construction. The main topic of this study is therefore the derivation of a tensor train representation of the first factor in each case, containing the action of the generator $\mathcal{L}$ (see Section~\ref{sec:tensor_reps}), and the efficient calculation of the reduced matrix by contraction of a suitable tensor network (see Section~\ref{sec:algorithms}).

\section{Tensor Representations of the Generator Data Tensor}
\label{sec:tensor_reps}
The goal of this section is to determine a tensor train representation of the generator data tensors $\mathcal{L}\PSI(x) \in \mathbb{R}^{n_1 \times \ldots \times n_p}$ and $\nabla \PSI(x) \in \mathbb{R}^{n_1 \times \ldots \times n_p \times d}$, given entry-wise by
\begin{align}
\label{eq:generator_data_tensors}
\mathcal{L}\PSI(x)_{u_1, \ldots, u_p} &= \mathcal{L}(\psi_{1, u_1} \cdot \ldots \cdot \psi_{p, u_p})(x), &
\nabla \PSI(x)_{u_1, \ldots, u_p, i} &= \pd{}{x^i} (\psi_{1, u_1} \cdot \ldots \cdot \psi_{p, u_p})(x).
\end{align}

\subsection{A General Representation Formula}
\label{subsec:general_representation}
As the Kolmogorov operator acts as a differential operator, its point-wise action on a product basis will be a sum of rank-one tensors, by the product rule. In anticipation of later results, we prove the following general representation formula, which is a straightforward extension of results shown in~\cite{Kazeev2013}:

\begin{proposition}
\label{prop:general_tt_decomposition}
Consider a tensor $\T \in \R^{n_1\times\dots\times n_p}$ defined by
		\begin{align}
		\label{eq:tensor_sum_rankone}
			\T =& \sum_{k=1}^{p}\big( \e_1 \otimes \dots \otimes \f_k \otimes \dots \otimes \e_p \big)\\
			\nonumber &+ \sum_{k_1=1}^{p-1} \sum_{k_2=(k_1 + 1)}^{p} \sum_{i=1}^{d} \big( \e_1 \otimes \dots \otimes \g_{k_1,i} \otimes \dots \otimes \h_{k_2,i} \otimes \dots \otimes \e_p  \big),
\end{align}
with vectors  $\e_k, \, \f_k, \, \g_{k,i}, \,\h_{k,i} \in \R^{n_k}$ for $k=1,\dots,p$ and $i=1,\dots,d$.
Then $\T$ can be written as a	tensor train $\T = \core{\T\super{1}}  \otimes \dots \otimes \core{\T\super{p}}$ with all ranks equal to $d + 2$, and cores:
		\begin{align}
		\label{eq:cores_general_tt_1}
			\core{\T\super{1}} &= \core{\begin{matrix}
					\e_1 & \f_1 & \g_{1,1} & \dots & \g_{1,d}
			\end{matrix}} \\
			\label{eq:cores_general_tt_2}
			\core{\T\super{k}} &= \core{\begin{matrix}
					\e_k & \f_k   & \g_{k, 1} & \dots & \g_{k,d} \\
					     & \e_k     &           &              &               \\
					     & \h_{k,1} & \e_k &              &     	        \\
					     & \vdots  &         &  \ddots  &              \\
					     & \h_{k,d}&        &               & \e_k 					     
			\end{matrix}} \qquad k=2,\dots,(p-1)\\
			\label{eq:cores_general_tt_3}
		\core{\T\super{p}} &= \core{\begin{matrix}
				\f_p &
				\e_p &
				\h_{p,1} &
				\dots &
				\h_{p,d}
		\end{matrix}}^T.
		\end{align}
\end{proposition}

\begin{proof}
	Define auxiliary tensors $\A\super{k},\, \B\super{k}_i \in \R^{n_k \times \dots \times n_p}$ for $k=1,\dots, (p-1)$ and $i=1,\dots,d$ by:
	\begin{align*}
		\A\super{k}=& \sum_{l=k}^{p}\big( \e_k \otimes \dots \otimes \f_l \otimes \dots \otimes \e_p \big)\\
			&+ \sum_{l_1=k}^{p-1} \sum_{l_2=(l_1 + 1)}^{p} \sum_{i=1}^{d} \big( \e_k \otimes \dots \otimes \g_{l_1,i} \otimes \dots \otimes \h_{l_2,i} \otimes \dots \otimes \e_p  \big),
	\end{align*}
	and
	\begin{align*}
		\B\super{k}_i &= \sum_{l=k}^{p}\big( \e_k \otimes \dots \otimes \h_{l,i} \otimes \dots \otimes \e_p \big).
	\end{align*}
	We also define $\A\super{p} = \f_p$ and $\B\super{p}_i = \h_{p, i}$. Note that $\A\super{k}$ is similar to $\T$, but all sums start from $k$ instead of $1$. In particular, we have that $\A\super{1} = \T$. The following relation holds for $k=1, \dots, (p-1)$, and is central for the proof:
	\begin{align} \label{eq:recursion_auxiliary}
		\A\super{k} = \mathbf{e}_k \otimes \A\super{k + 1} + \big( \f_k \otimes \e_{k+1} \otimes \dots \otimes \e_p \big) + \sum_{i=1}^d \g_{k,i} \otimes \B\super{k+1}_i.
	\end{align}
The recursion~\eqref{eq:recursion_auxiliary} can be verified directly by expanding $\mathbf{e}_k \otimes \A\super{k + 1}$, and noting that the remaining terms in~\eqref{eq:recursion_auxiliary} are exactly what is missing to complete the expression for $\A\super{k}$. We also obtain a recursive relation for $\B\super{k}_i$, for $k = 1, \ldots, (p-1)$:

\begin{align}
\label{eq:recursion_Bk}
\B\super{k}_i &= \h_{k, i} \otimes \e_{k+1} \otimes \dots \otimes \e_p + \e_k \otimes \B\super{k+1}_i.
\end{align}
Using (\ref{eq:recursion_auxiliary} - \ref{eq:recursion_Bk}), we can verify that for $k = 2, \dots, p-1$:
	\begin{align} &\nonumber
		\core{\T\super{k}} \otimes \core{\begin{matrix}
				\A\super{k+1} &
				\e_{k+1} \otimes \dots \otimes \e_p &
				\B\super{k+1}_1 &
				\dots &
				\B\super{k+1}_d
			\end{matrix}}^T \\
		&= \core{\begin{matrix}
					\e_k & \f_k   & \g_{k, 1} & \dots & \g_{k,d}  \\
					     & \e_k     &           &              &               \\
					     & \h_{k,1} & \e_k &              &     	       \\
					     & \vdots  &         &  \ddots  &               \\
					     & \h_{k,d}&        &               & \e_k      
			\end{matrix}}
		\otimes \core{\begin{matrix}
				\A\super{k+1} \\
				\e_{k+1} \otimes \dots \otimes \e_p\\
				\B\super{k+1}_1\\
				\vdots\\
				\B\super{k+1}_d
			\end{matrix}} \\
	&= \core{\begin{matrix}
			\A\super{k} &
			\e_{k} \otimes \dots \otimes \e_p &
			\B\super{k}_1 &
			\dots &
			\B\super{k}_d
	\end{matrix}}^T. \label{eq:recursion_tensordor}
	\end{align}
By noting that 
\[ \core{\T\super{p}} = \core{\begin{matrix}
			\A\super{p} &
			 \e_p &
			\B\super{p}_1 &
			\dots &
			\B\super{p}_d
	\end{matrix}}^T, \]
and repeated application of~\eqref{eq:recursion_tensordor}, we find
	\begin{align*}
		\core{\T\super{1}}  \otimes \dots \otimes \core{\T\super{p}} 
		&= \core{\begin{matrix}
				\e_1 & \f_1 & \g_{1,1} & \dots & \g_{1,d}
		\end{matrix}} \otimes \core{\begin{matrix}
				\A\super{2} \\
				\e_{2} \otimes \dots \otimes \e_p \\
				\B\super{2}_1 \\
				\vdots\\
				\B\super{2}_d
		\end{matrix}} 
	= \T.
	\end{align*}
\end{proof}

\begin{remark}
Proposition~\ref{prop:general_tt_decomposition} differs from the setting in~\cite{Kazeev2013} by adding the sums over $i$ in the second line of~\eqref{eq:tensor_sum_rankone}. On the other hand, additional weights $a_{k_1, k_2}$ in each term were considered in~\cite{Kazeev2013}, as well as an additional sum over the terms in the second line, but with $\g_{k, i}$ and $\h_{k, i}$ swapped. The latter modification, without weights, is easily incorporated into~(\ref{eq:cores_general_tt_1} - \ref{eq:cores_general_tt_3}). However, this result is not needed in the following:
\end{remark}

\begin{corollary}
\label{cor:symmetric_tt_rep}
If the definition of $\T$ additionally contains the terms:
\[ \sum_{k_1=1}^{p-1} \sum_{k_2=(k_1 + 1)}^{p} \sum_{i=1}^{d} \big( \e_1 \otimes \dots \otimes \h_{k_1,i} \otimes \dots \otimes \g_{k_2,i} \otimes \dots \otimes \e_p  \big), \]
then $\T$ can be represented with cores of ranks $2d + 2$ similar to~(\ref{eq:cores_general_tt_1} - \ref{eq:cores_general_tt_3}). In detail, $\T\super{1}$ is extended by appending all terms $\h_{1, i}$, while all $\g_{p, i}$ are appended to $\T\super{p}$. For $2 \leq k \leq p-1$, all terms $\h_{k, i}$ are added to the first row of $\T\super{k}$, while all $\g_{k, i}$ are added to the second column. The diagonal is just filled up with copies of the vector $\e_{k}$.
\end{corollary}

\subsection{Representations of the Generator Data Tensor}
\label{subsec:rep_op_data_tensor}
Based on Proposition~\ref{prop:general_tt_decomposition}, we can now derive tensor train representations of the data tensors required for generator EDMD in TT format. We introduce the following notation: for a vector-valued function $\psi : \mathbb{X} \mapsto \mathbb{R}^{n}$, we write $\frac{\partial \psi}{\partial x^i}$ for the vector of partial derivatives of all functions $\psi_{u}$ with respect to coordinate $x^i$. The associated Jacobian matrix is denoted by $\nabla \psi = \left(\frac{\partial \psi_{u}}{\partial x^i}\right)_{u, i} \in \mathbb{R}^{n \times d}$, and the corresponding Hessian tensor by $\nabla^2 \psi \in \mathbb{R}^{n \times d \times d}$.

\subsubsection{Reversible Case}
\label{subsubsec:rep_op_rev}
We consider the reversible generator $\mathcal{L}$ first, since only first order derivatives are required for its data-driven approximation \eqref{eq:reduced_matrix_pw_rev}.

\begin{lemma}
\label{lem:tt_decomp_dpsi_x}
For any $x \in \mathbb{X}$, the tensor $\nabla \PSI(x) \in \mathbb{R}^{n_1 \times \ldots \times n_p \times d}$ can be written in the form~\eqref{eq:tensor_sum_rankone} with:
\begin{align*}
\e_k &= \begin{cases} \p_k(x) & 1 \leq k \leq p \\ 0 & k = p+1 \end{cases},
& \f_k &= 0, & \g_{k, i} &= \pd{\p_k}{x^i}(x), &
\h_{k, i}  &= \begin{cases} 0 & 2\leq k \leq p \\ e_i & k=p+1. \end{cases}
\end{align*}
\end{lemma}
\begin{proof}
The representation follows directly from the product rule:
\begin{align*}
\nabla \PSI(x) &= \sum_{i=1}^d \sum_{k = 1}^p \p_1(x) \otimes \dots \otimes \pd{\p_k}{x^i}(x) \otimes \dots \otimes \p_p(x) \otimes e_i,
\end{align*}
which is already~\eqref{eq:tensor_sum_rankone} with the parameters given by the Lemma.
\end{proof}

\begin{corollary}
\label{cor:tt_nabla_psi}
The tensor train representation of $\nabla \PSI(x)$, with ranks equal to $r_1 = (d + 1), \ldots, r_{p-1} = (d+1), r_p = d$, is given by the cores:
\begin{align*}
\core{\nabla \PSI\super{1}(x)} &= \core{\begin{matrix}
					\p_1(x) & \pd{\psi_1}{x^1}(x) & \dots & \pd{\psi_1}{x^d}(x)
					\end{matrix}} \\
			\core{\nabla \PSI\super{k}(x)} &= \core{\begin{matrix}
						\p_k(x) & \pd{\psi_k}{x^1}(x) & \dots & \pd{\psi_k}{x^d}(x)  \\
					     & \p_k(x)     &           &        \\
					     &   				&    \ddots     &              \\
					     &  			&       			 &  \p_k(x) 
			\end{matrix}}, \quad 2 \leq k \leq p - 1, \\
			\core{\nabla \PSI\super{p}(x)} &= \core{\begin{matrix}
						\pd{\psi_p}{x^1}(x) & \dots & \pd{\psi_p}{x^d}(x)  \\
					     \p_p(x)     &           &        \\
					      				 &    \ddots     &              \\
					     			&       			 &  \p_p(x) 
			\end{matrix}}, \\
		\core{\nabla \PSI\super{p+1}(x)} &= \core{\begin{matrix}
				e_1 &
				\dots &
				e_d
		\end{matrix}}^T.
\end{align*}
\end{corollary}
\begin{proof}
The result follows from Lemma~\ref{lem:tt_decomp_dpsi_x}, and by observing that the second row and column of each core in Proposition~\ref{prop:general_tt_decomposition} can be omitted in this case, as they do not contribute to the final tensor. Moreover, the first column of the pre-last core $\nabla \PSI\super{p}(x)$ can also be omitted, as it would be multiplied by a zero in $\nabla \PSI\super{p+1}(x)$.
\end{proof}

\begin{remark}
A TT representation of $\mathrm{d}\PSI(x)$ can be obtained based on Corollary~\ref{cor:tt_nabla_psi} by absorbing the additional factor $\sigma$ into the final core.
As $\mathrm{d}\PSI(X) = \sum_{l=1}^m \mathrm{d}\PSI(x_l) \otimes e_l$, where $e_l$ is the canonical unit vector in $\mathbb{R}^m$, we can also obtain a TT decomposition of $\mathrm{d}\PSI(X)$ using the summation rule for tensor trains \cite{Oseledets2011}.
\end{remark}

\subsubsection{Non-reversible Case}
\label{subsubsec:rep_op_non_rev}
Next, we consider the generator $\mathcal{L}$ for a general (not necessarily reversible) SDE. The point-wise evaluation of $\mathcal{L}$ can be written as a sum of rank-one tensors as follows:
\begin{lemma}
\label{lem:tt_decomp_Lpsi_x}
For any $x \in \mathbb{X}$, the tensor $\mathcal{L}\PSI(x) \in \mathbb{R}^{n_1 \times \ldots \times n_p}$ can be written in the form \eqref{eq:tensor_sum_rankone} with:
\begin{align*}
\e_k &= \p_k(x), & \f_k &= \mathcal{L}\psi_k(x), & \g_{k, i} &= \h_{k, i}  = \nabla \psi_k(x) \cdot \sigma_{:, i}(x),
\end{align*}
where $\sigma_{:, i}$ denotes the $i$-th column of the diffusion $\sigma$.
\end{lemma}

\begin{proof}
The first and second derivative of the tensor $\PSI(x)$ are given by (omitting the argument $x$):
\begin{align*}
\pd{\PSI}{x^i} &= \sum_{k=1}^p \p_1 \otimes \dots \otimes \pd{\psi_k}{x^i} \otimes \dots \otimes \p_p, \\
\pdtwo{\PSI}{x^i}{x^j} &= \sum_{k=1}^p \p_1 \otimes \dots \otimes \pdtwo{\psi_k}{x^i}{x^j} \otimes \dots \otimes \p_p + \sum_{k_1 \neq k_2} \p_1 \otimes \dots \otimes \pd{\psi_{k_1}}{x^i} \otimes \dots \otimes \pd{\psi_{k_2}}{x^j} \otimes \dots \otimes \p_p.
\end{align*}
Therefore, application of the generator yields:
\begin{align*}
\mathcal{L}\PSI &= \sum_{i=1}^d b_i \pd{\PSI}{x^i} + \frac{1}{2} \sum_{i,j=1}^d a_{ij} \pdtwo{\PSI}{x^i}{x^j} \\
&= \sum_{k=1}^p \p_1 \otimes \dots \otimes (\nabla \psi_k \cdot b + \frac{1}{2} \nabla^2 \psi_k : a)  \otimes \dots \otimes \p_p  \\
&+ \frac{1}{2}\sum_{i,j=1}^d a_{ij} \sum_{k_1 \neq k_2} \p_1 \otimes \dots \otimes \pd{\psi_{k_1}}{x^i} \otimes \dots \otimes \pd{\psi_{k_2}}{x^j} \otimes \dots \otimes \p_p.
\end{align*}
The $k$-th term in the first line is just $\mathcal{L}\psi_k$, while in the second line, we can write $a_{ij} = \sum_{s=1}^d \sigma_{is}\sigma_{js}$, to obtain:
\begin{align*}
\mathcal{L}\PSI &= \sum_{k=1}^p \p_1 \otimes \dots \otimes \mathcal{L} \psi_k  \otimes \dots \otimes \p_p  \\
&+ \frac{1}{2} \sum_{s=1}^d \sum_{k_1 \neq k_2} \p_1 \otimes \dots \otimes (\nabla \psi_{k_1} \cdot \sigma_{:, s}) \otimes \dots \otimes (\nabla \psi_{k_2} \cdot \sigma_{:, s}) \otimes \dots \otimes \p_p \\
&= \sum_{k=1}^p \p_1 \otimes \dots \otimes \mathcal{L} \psi_k  \otimes \dots \otimes \p_p  \\
&+ \sum_{s=1}^d \sum_{k_1=1}^{p-1} \sum_{k_2=k_1+1}^p \p_1 \otimes \dots \otimes (\nabla \psi_{k_1} \cdot \sigma_{:, s}) \otimes \dots \otimes (\nabla \psi_{k_2} \cdot \sigma_{:, s}) \otimes \dots \otimes \p_p.
\end{align*}
In the last step, we used that the second term is symmetric in $k_1$ and $k_2$. This is the desired result.
\end{proof}
\noindent The tensor train representation of $\mathcal{L}\PSI(x)$ is then obtained from Proposition~\ref{prop:general_tt_decomposition}:

\begin{corollary}
\label{cor:tt_rep_Lpsi_non_rev}
The tensor train representation of $\mathcal{L}\PSI(x)$, with ranks all equal to $(d+2)$, is given by the cores:
\begin{align*}
\core{\mathcal{L}\PSI\super{1}(x)} &= \core{\begin{matrix}
					\p_1(x) & \mathcal{L}\psi_1(x) & \nabla \psi_1(x) \cdot \sigma_{:, 1}(x) & \dots & \nabla \psi_1(x) \cdot \sigma_{:, d}(x) \end{matrix}} \\
			\core{\mathcal{L}\PSI\super{k}(x)} &= \core{\begin{matrix}
						\p_k(x) & \mathcal{L}\psi_k(x) & \nabla \psi_k(x) \cdot \sigma_{:, 1}(x) & \dots & \nabla \psi_k(x) \cdot \sigma_{:, d}(x) \\
					     & \p_k(x)     &           &              &     \\
					     & \nabla \psi_k(x) \cdot \sigma_{:, 1}(x) & \p_k(x) &              &      \\
					     & \vdots  &         &  \ddots  &            \\
					     & \nabla \psi_k(x) \cdot \sigma_{:, d}(x) &        &               & \p_k(x) 
			\end{matrix}} \\
		\core{\mathcal{L}\PSI\super{p}(x)} &= \core{\begin{matrix}
				\mathcal{L}\p_p(x) &
				\p_p(x) &
				\nabla \psi_p(x) \cdot \sigma_{:, 1}(x) &
				\dots &
				\nabla \psi_p(x) \cdot \sigma_{:, d}(x)
		\end{matrix}}^T.
\end{align*}
\end{corollary}

\begin{remark}
As $\mathcal{L}\PSI(X) = \sum_{l=1}^m \mathcal{L}\PSI(x_l) \otimes e_l$, where $e_l$ is the canonical unit vector in $\mathbb{R}^m$, we can also obtain a TT decomposition of $\mathcal{L}\PSI(X)$ using the summation rule for tensor trains \cite{Oseledets2011}.
\end{remark}

\section{Algorithmic Realization}
\label{sec:algorithms}
In this section, we discuss the algorithmic realization of the tensor-based version of gEDMD. We first present the tgEDMD algorithm in Section~\ref{subsec:efficient_calculation_M} and briefly analyze its complexity and consistency for infinite data. Subsequently, we present two extensions of the method in Section~\ref{subsec:extensions}.

\subsection{The tgEDMD Algorithm}
\label{subsec:efficient_calculation_M}
Based on the results of Section~\ref{sec:tensor_reps}, the tensor-based counterpart of the AMUSE Algorithm~\ref{alg:AMUSE}, which we call tgEDMD, is given in Algorithm~\ref{alg:tgedmd} below. Its output is a matrix representation $\hat{M}_{r_p}$ or $\hat{M}_{r_p}^{\mathrm{rev}}$ of the generator, on the linear span of the functions $\hat{\eta}_{r_p}$. This matrix representation can then be used for further analysis of the generator, e.g. for spectral analysis, prediction, or model reduction.

\begin{algorithm}[htb]
  \caption{tgEDMD}
  \label{alg:tgedmd}
  \setlength{\tabcolsep}{.5ex}
  \begin{tabular}{ll}
    \textbf{Input:} & Data $x_1, \ldots, x_m \in \mathbb{X}$, basis sets $\{\psi_{k, u_k} \}_{u_k=1}^{n_k}, \, 1\leq k \leq p$. \\
    \textbf{Output:} & Reduced matrix representation $\hat{M}_{r_p}\in \mathbb{R}^{r_p \times r_p}$ or $\hat{M}_{r_p}^{\mathrm{rev}} \in \mathbb{R}^{r_p \times r_p}$ of $\mathcal{L}$.
  \end{tabular}
  \hrule\vspace{0.2cm}
  \begin{algorithmic}[1]
    \State Compute global SVD of $\PSI(X)$, i.e. $\PSI(X)\vert_p \approx \hat{\mathbf{U}}\vert_p \hat{\Sigma} \hat{V}^T$, where $\hat{\mathbf{U}}$ is in tensor train format with ranks $r_1, \ldots, r_p$ and left-orthonormal cores.
    \State Set up TT representations of $\mathcal{L}\PSI(x_l)$ or $\nabla \PSI(x_l)$ for $1 \leq l \leq m$, according to Corollaries~\ref{cor:tt_rep_Lpsi_non_rev} or~\ref{cor:tt_nabla_psi}.
    \State Compute reduced matrices $\hat{M}_{r_p}$ or $\hat{M}_{r_p}^{\mathrm{rev}}$, according to~(\ref{eq:amuse_approx_generator} - \ref{eq:amuse_approx_generator_rev}).
  \end{algorithmic}
\end{algorithm}

\paragraph{Computational Effort} The computational effort required by Algorithm~\ref{alg:tgedmd} can be broken down into three major contributions.

\noindent The first is the global SVD of the data tensor $\PSI(X)$. The method consists of a series of matrix SVDs, each matrix being of size $r_{k-1} n_k \times m$. As is well known, the cost  for these can be estimated as $\mathcal{O}(\min\{r_{k-1}^2 n_k^2 m, r_{k-1} n_k m^2\})$. The cost of updating the next core in each step can be neglected in comparison.

\noindent Secondly, the reduced matrices $\hat{M}_{r_p}$ or $\hat{M}_{r_p}^{\mathrm{rev}}$ need to be computed. To this end, for each data point $x_l$, the contraction of the orthonormal tensor train $\hat{\mathbf{U}}$ with either $\mathcal{L}\PSI(x_l)$ or $\nabla \PSI(x_l)$ in~(\ref{eq:amuse_approx_generator} - \ref{eq:amuse_approx_generator_rev}) must be carried out. This can be achieved by Algorithm~\ref{alg:contraction_tt_network} below, see~\cite{Oseledets2011} for details. The tensor product in lines~\ref{algline:kronecker_1} and~\ref{algline:kronecker_k} is the Kronecker product between matrices defined in~\eqref{eq:kronecker_product}. Of course, Algorithm~\ref{alg:contraction_tt_network} can be applied to any pair of general tensors of compatible dimensions.

\begin{algorithm}[htb]
  \caption{Tensor Network Contraction}
  \label{alg:contraction_tt_network}
  \setlength{\tabcolsep}{.5ex}
  \begin{tabular}{ll}
    \textbf{Input:} & Tensor trains $\mathbf{T}, \hat{\mathbf{U}}$ of order $p$ with ranks $s_1, \ldots, s_p$ and $r_1, \ldots, r_p$, respectively, \\
    	& where either $\mathbf{T} = \mathcal{L}\PSI(x_l)$ or $\mathbf{T} = \nabla \PSI(x_l)$. \\
    \textbf{Output:} & Contraction $\mathbf{T} \cdot \hat{\mathbf{U}} \in \mathbb{R}^{s_p \times r_p}$, defined by $(\mathbf{T} \cdot \hat{\mathbf{U}})_{u, v} = \sum_{u_1,\ldots, u_p} \mathbf{T}_{u_1, \ldots, u_p; u} \, \hat{\mathbf{U}}_{u_1, \ldots, u_p; v}$.
  \end{tabular}
  \hrule\vspace{0.2cm}
  \begin{algorithmic}[1]
	\State $ v = \sum_{u_1} \mathbf{T}\super{1}_{:, u_1, :} \otimes \hat{\mathbf{U}}\super{1}_{:, u_1, :}$.	  \label{algline:kronecker_1}	
	\For{$k = 2, \ldots, p$}
		\State $ v \gets v \left[\sum_{u_k} \mathbf{T}\super{k}_{:, u_k, :} \otimes \hat{\mathbf{U}}\super{k}_{:, u_k, :}\right]$.
		\label{algline:kronecker_k}
	\EndFor
	\State Re-shape $v$ to shape $s_p \times r_p$.
	\State \Return $v$
  \end{algorithmic}
\end{algorithm}

\noindent As the cores of $\mathcal{L}\PSI(x_l)$ and $\nabla \PSI(x_l)$ are structured and sparse, we obtain the following result:

\begin{lemma}
The cost of applying Algorithm~\ref{alg:contraction_tt_network} can be estimated by $\mathcal{O}(\sum_{k=1}^p n_k d r_{k-1}r_k)$.
\end{lemma}
\begin{proof}
We obtain the result by recalling that the TT ranks of $\mathcal{L}\PSI(x_l)$ and $\nabla \PSI(x_l)$ are essentially equal to $d$ by the results obtained in Sec.~\ref{sec:tensor_reps}, and by observing that the sparsity patterns of their cores are invariant under Kronecker products. Taking $\mathcal{L}\PSI(x_l)$ as an example, we have by~\eqref{eq:kronecker_product} for all $2 \leq k \leq p-1$:
\begin{align*}
\mathcal{L}\PSI(x_l)\super{k}_{:, u_k, :} \otimes \hat{\mathbf{U}}\super{k}_{:, u_k, :} = \begin{bmatrix}
B_{1,1} & B_{1,2} & B_{1,3} & \dots & B_{1, d+2} \\
			& B_{2,2} & 			 &			  &				 \\
			& B_{3,2} & B_{3,3} & 			&				 \\
			&  \vdots & 			&  \ddots &				\\
			& B_{d+2, 2} & 		&				& B_{d+2, d+2}
\end{bmatrix},
\end{align*}
where each block $B_{u, v} \in \mathbb{R}^{r_{k-1} \times r_k}$ is given by $\mathcal{L}\PSI(x_l)\super{k}_{u, u_k, v} \, \hat{\mathbf{U}}\super{k}_{:, u_k, :}$. The Kronecker product $\mathcal{L}\PSI(x_l)\super{k}_{:, u_k, :} \otimes \hat{\mathbf{U}}\super{k}_{:, u_k, :}$ can therefore be evaluated in $\mathcal{O}(d r_{k-1} r_k)$ operations, while multiplication by a row vector from the left also costs $\mathcal{O}(d r_{k-1} r_k)$ operations. The full cost of applying Algorithm~\ref{alg:contraction_tt_network} is therefore $\mathcal{O}(\sum_{k=1}^p n_k d r_{k-1}r_k)$.
\end{proof}

\noindent The remaining matrix multiplications appearing in the computation of the reduced matrices are again not relevant in comparison, leaving us with the estimate $\mathcal{O}(m\sum_{k=1}^p n_k d r_{k-1}r_k)$. Finally, if eigenvalues of the Kolmogorov operator are required, an additional cost of $\mathcal{O}(r_p^3)$ for the diagonalization of the reduced matrix must be allowed. In summary, the total cost for the approximation of eigenvalues of the Kolmogorov operator can be estimated as
\begin{equation}
\label{eq:comp_cost_tgedmd}
\mathcal{O}(\sum_{k=1}^p \min\{r_{k-1}^2 n_k^2 m, r_{k-1} n_k m^2\}) + \mathcal{O}(m\sum_{k=1}^p n_k d r_{k-1}r_k) + \mathcal{O}(r_p^3).
\end{equation}

\noindent It is instructive to compare these findings to the cost of applying the standard AMUSE Algorithm~\ref{alg:AMUSE} in conjunction with the full product basis $\PSI$. Denoting the product of all mode sizes by $N = \prod_{k=1}^p n_k$, and adding up the operation counts for the SVD of the data tensor $\PSI(X)$, the assembly of the reduced matrix, and diagonalization of the latter, we end up with the estimate
\begin{equation}
\label{eq:comp_cost_gedmd}
\mathcal{O}(N m^2) + \mathcal{O}(m(N r + r^2)) + \mathcal{O}(r^3).
\end{equation}
In all but the smallest examples, the cost is dominated by the exponential dependence on the dimension $p$ through $N$.

\paragraph{Consistency}
The discussion on consistency of the AMUSE algorithm, see Section~\ref{subsec:AMUSE}, carries over to the tensor case, based on the analysis in Ref.~\cite[Section 5]{Nueske2021}. Therein, it was shown that if all SVDs in the global SVD are truncated according to a fixed vector of TT ranks $\mathbf{r}  = [r_1, \ldots, r_p]$, one can construct a recursive sequence of subspaces $\mathbb{G}^{:k}, \, 1\leq k \leq p$ (called multi-linear spectral subspaces), such that the linear span of the functions $\hat{\eta}_{r_p}$ consistently approximates the last of these subspaces $\mathbb{G}^{:p}$, along with the projected Koopman generator. This result required non-degeneracy of a series of eigenvalues of recursively constructed Gramian matrices. Analogous to the matrix case, we can also impose a vector of truncation parameters $[\epsilon_1, \ldots, \epsilon_p]$ instead and choose the reduced rank of each SVD in the global SVD method such that $\hat{\sigma}^k_{r_k} \geq \sqrt{m} \epsilon_k > \hat{\sigma}^k_{r_k + 1}, \, 1\leq k \leq p$. Given similar non-degeneracy conditions, the same arguments as in~\cite{Nueske2021} can be used to show that there are limiting ranks $r_k(\epsilon_k)$ such that $r_k  \rightarrow r_k(\epsilon_k), \, 1\leq k \leq p$, and the linear span of $\hat{\eta}_{r_p}$ consistently approximates $\mathbb{G}^{:p}$ for the TT ranks $[r_1(\epsilon_1), \ldots, r_p(\epsilon_p)]$.

\subsection{Extensions}
\label{subsec:extensions}
We conclude the methodological part of this paper by briefly discussing two practically relevant settings which require modifications of the proposed methods.

\paragraph{Importance Sampling}
First, obtaining data sampling the invariant measure $\mu$ is often a major challenge, for instance due to metastability of the dynamics \eqref{eq:sde}. In principle, this difficulty can be circumvented  as long as appropriate re-weighting factors can be calculated. Assume the data are generated according to a probability distribution with density $\vartheta$ on $\mathbb{X}$ and we are able to calculate the \emph{importance sampling ratio} $w(x) = Z\frac{\rho(x)}{\vartheta(x)}$ at any point $x\in \mathbb{X}$, where $Z$ is a possibly unknown constant. This setting reflects a rather typical scenario where both the invariant density $\rho$ and the sampling density $\vartheta$ can be evaluated up to an unknown normalization constant. Note, however, that a poor choice of the sampling density $\vartheta$ can drastically increase the variance of all estimators involved, but this problem is beyond the scope of the present work.

\noindent To see how this affects the proposed methods, consider the matrix case first. Introducing the diagonal weighting matrix $W = \mathrm{diag}[w(x_1), \ldots, w(x_m)]$, the empirical estimators (\ref{eq:empirical_estimator_C} - \ref{eq:empirical_estimator_A}) change according to:
\begin{align*}
\hat{C}(\psi) &= \frac{1}{m}(\Psi(X) W^{1/2}) (\Psi(X) W^{1/2})^T, & \hat{A}(\psi) &= \frac{1}{m}(\Psi(X) W^{1/2}) (\mathcal{L}\Psi(X) W^{1/2})^T.
\end{align*}
We need to modify line~\ref{algline:AMUSE_SVD} of Algorithm~\ref{alg:AMUSE} to the computation of the rank-$r$ SVD $(\Psi(X)W^{1/2}) \approx \hat{U}_r \hat{\Sigma}_r \hat{V}_r^T$. The reduced matrices are then obtained as
\begin{align*}
\hat{M}_r &= \hat{\Sigma}_r^{-1} \, \hat{U}_r^\top \, (\Psi(X) W^{1/2}) (\, \mathcal{L}\Psi(X)W^{1/2})^T \, \hat{U}_r \, \hat{\Sigma}_r^{-1} = \hat{V}_r^\top W^{1/2} \, \mathcal{L}\Psi(X)^T \, \hat{U}_r \, \hat{\Sigma}_r^{-1}, \\
\hat{M}_r^{\mathrm{rev}} &= -\frac{1}{2} \hat{\Sigma}_r^{-1} \, \hat{U}_r^\top \, \mathrm{d}\Psi(X) \, W \, \mathrm{d}\Psi(X)^T  \, \hat{U}_r \, \hat{\Sigma}_r^{-1}.
\end{align*}
Note that both estimators are indeed invariant with respect to the constant $Z$. In the tensor case, the following modifications must be made as a result:
\begin{enumerate}
\item The TT representation~\eqref{eq:tt_decomp_psi_x} of $\PSI(X)$ is modified by scaling each unit vector in the last core by the corresponding re-weighting factor:
\[ \core{\PSI(X)\super{p+1}}  = \core{\begin{matrix}
\sqrt{w(x_1)} e_1 \\ \vdots \\ \sqrt{w(x_m)} e_m
\end{matrix}}. \]
\item Application of the global SVD of $\PSI(X)$ must be preceded by right-orthonormalization of $\PSI(X)\super{p+1}$, see \cite{KLUS2018}.
\item Reduced matrices are calculated according to
\begin{align*}
\hat{M}_{r_p} &= \sum_{l=1}^m w(x_l)^{1/2} \hat{V}_{l, :} \otimes (\mathcal{L}\PSI(x_l)\vert_p)^T \hat{\mathbf{U}}\vert_p \hat{\Sigma}^{-1}, \\
\nonumber
\hat{M}_{r_p}^{\mathrm{rev}} &= -\frac{1}{2} \sum_{l=1}^m w(x_l) \hat{\Sigma}^{-1} \, (\hat{\mathbf{U}}\vert_p)^T \, \nabla \PSI(x_l) \,a(x_l) \, \nabla \PSI(x_l)^T \, \hat{\mathbf{U}}\vert_p \, \hat{\Sigma}^{-1}.
\end{align*}
\end{enumerate} 

\paragraph{Projection onto Reduced Variables}
Secondly, we discuss the use of basis functions defined on a set of reduced variables. For high-dimensional systems, it is often convenient to employ basis functions which are not directly defined on the state variables $x^i, \, 1 \leq i \leq d$, but on the image of some (possibly) non-linear descriptors 
\[ \xi: \mathbb{X} \mapsto \mathbb{Y} \subset \mathbb{R}^D, \, D \leq d. \] 
Now, if $\mathbb{V}  =\mathrm{span}\{\psi_u(y) \}_{u=1}^n$ is a finite-dimensional space of functions defined on the lower-dimensional state space $\mathbb{Y}$, the resulting Galerkin approximation to the generator possesses a useful interpretation, as discussed in great detail in Ref.~\cite{Zhang2016}. Consider the space of functions in $L^2_\mu$ which are in fact functions of $\xi$, i.e.
\[ \mathbb{V}^\xi = \{\bar{\phi} \in L^2_\mu(\mathbb{X}): \quad \exists \phi: \mathbb{Y} \mapsto \mathbb{R}, \, \bar{\phi} = \phi \circ \xi \}.\]
After identification of $\bar{\phi}$ with $\phi$, $\mathbb{V}^\xi$ equals the weighted Lebesgue space $L^2_\nu(\mathbb{Y})$, where $\nu$ is the marginal probability distribution of $\mu$ along $\xi$. Using the $\mu$-orthogonal projector $\mathcal{P}^\xi$ onto $L^2_\nu(\mathbb{Y})$, one can introduce the projected generator $\mathcal{L}^\xi = \mathcal{P}^\xi \mathcal{L}\mathcal{P}^\xi$. Given some technical assumptions \cite{Zhang2016}, the projected generator $\mathcal{L}^\xi$ is again the generator of a stochastic differential equation \eqref{eq:sde}, with effective drift and diffusion coefficients. Importantly, the following identities relating the generators $\mathcal{L}, \, \mathcal{L}^\xi$ and their invariant measures $\mu, \, \nu$, hold for all $\phi, \tilde{\phi} \in L^2_\nu(\mathbb{Y})$:
\begin{align*}
\innerprod{\phi}{\tilde{\phi}}_\nu &= \innerprod{\phi\circ \xi}{\tilde{\phi}\circ \xi}_\mu, & \innerprod{\phi}{\mathcal{L}^\xi \tilde{\phi}}_\nu &= \innerprod{\phi\circ \xi}{\mathcal{L}(\tilde{\phi}\circ \xi)}_\mu.
\end{align*}
\noindent As a consequence, any Galerkin approximation to $\mathcal{L}$ on a subspace $\mathbb{V} \subset L^2_\nu(\mathbb{Y})$ is simultaneously an approximation to $\mathcal{L}^\xi$. Conversely, the empirical estimators~(\ref{eq:empirical_estimator_C} - \ref{eq:empirical_estimator_Arev}) can be used in conjunction with simulation data of the full process, sampling $\mu$, in order to approximate the reduced generator $\mathcal{L}^\xi$ \cite{KLUS2020}. When applying the full generator $\mathcal{L}$ to a basis set on $\mathbb{Y}$ as above, derivatives of $\xi$ with respect to the full state variables $x^i$ must be taken using the chain rule.

\noindent If $\mathbb{V} \subset \mathcal{D}(\mathcal{L}^\xi)$ is a tensor subspace, the following adjustments need to be made when applying tgEDMD:
\begin{enumerate}
\item In the reversible case, it suffices to replace $a(x_l)$ in ~\eqref{eq:amuse_approx_generator_rev} by $a^\xi(x_l) = \nabla \xi(x_l) a(x_l) \nabla \xi(x_l)^T$. All derivatives appearing in Corollary~\ref{cor:tt_nabla_psi} can then be taken simply with respect to the reduced variables $y$.
\item In the non-reversible setting, the terms $\nabla \psi_k(x) \cdot \sigma_{:, i}(x)$ in Corollary~\ref{cor:tt_rep_Lpsi_non_rev} must be replaced by $\nabla_y \psi_k(y) \cdot \nabla \xi(x) \cdot \sigma_{:, i}(x)$, where $y = \xi(x)$. Moreover, each of the terms $\mathcal{L}\psi_k(x)$ must be evaluated as
\[ \mathcal{L}\psi_k(x) = \nabla_y \psi_k(y) \cdot \nabla \xi(x) \cdot b(x) + \frac{1}{2} \nabla^2_y \psi_k(y) : a^\xi(x) + \frac{1}{2} \nabla_y \psi_k(y) (\nabla^2 \xi(x) : a(x)). \]
\end{enumerate}

\section{Numerical Results}
\label{sec:numerical_results}
Below, we illustrate the capabilities of the proposed methods by means of two examples. The first is a diffusion in a four-dimensional model potential (Section~\ref{subsec:lemon_slice}), the second is a data set of molecular dynamics simulation of the deca-alanine peptide (Section~\ref{subsec:deca_alanine}). For both examples, we mainly focus on re-producing the dominant spectrum of the generator. In the setting outlined in Sec.~\ref{subsec:kolmogorov_op}, the generator $\mathcal{L}$ is non-negative, hence any eigenvalue $\kappa_k$ must satisfy $\mathcal{R}(\kappa_k) \leq 0$, where $\mathcal{R}(\cdot)$ is the real part. If the system is reversible in addition, all eigenvalues lie on the negative real half-line. For a broad class of systems, the largest eigenvalue $\kappa_0 = 0$ is non-degenerate, possibly followed by a number of real eigenvalues close to zero and a subsequent spectral gap, that is
\[ 0 = \kappa_0 > \kappa_1 \geq \ldots \geq \kappa_c \gg R, \]
where $R < 0$ is an upper bound for the remaining spectrum. This pattern, which we will encounter in both examples, is the signature of metastability, i.e. the existence of long-lived states such that interstate transitions are rare-events \cite{Davies1982}.

\noindent  By the spectral mapping theorem \cite{Pazy2012}, any generator eigenvalue $\kappa_k$ corresponds to an eigenvalue $\lambda_k(t) = e^{\kappa_k t}$ of the Koopman operator. If $\mathcal{R}(\kappa_k) < 0$, the Koopman eigenvalue is thus decaying exponentially with the time lag $t$. The characteristic decay time scale associated to $\lambda_k(t)$, also called \emph{implied time scale}, is
\[ t_k := -\frac{1}{\mathcal{R}(\kappa_k)} = -\frac{t}{\log |\lambda_k(t)|}, \]
where $|\cdot |$ is the modulus of a complex number. Implied time scales have frequently been used as a metric in order to compare different approximate models for the Koopman operator or generator \cite{Prinz2011c,BowmanPandeNoe_MSMBook}.

\noindent Our software implementation of the tgEDMD algorithm closely follows the discussion in Section~\ref{subsec:efficient_calculation_M}. It is publicly available as part of the scikit-tt package\footnote{\url{https://github.com/PGelss/scikit_tt/}}, including the numerical examples presented below. The data required to re-produce these examples is also available online\footnote{\url{https://doi.org/10.5281/zenodo.6367143}}.

\subsection{Lemon Slice Potential}
\label{subsec:lemon_slice}
\begin{figure}
\centering
\includegraphics[width=0.46\textwidth]{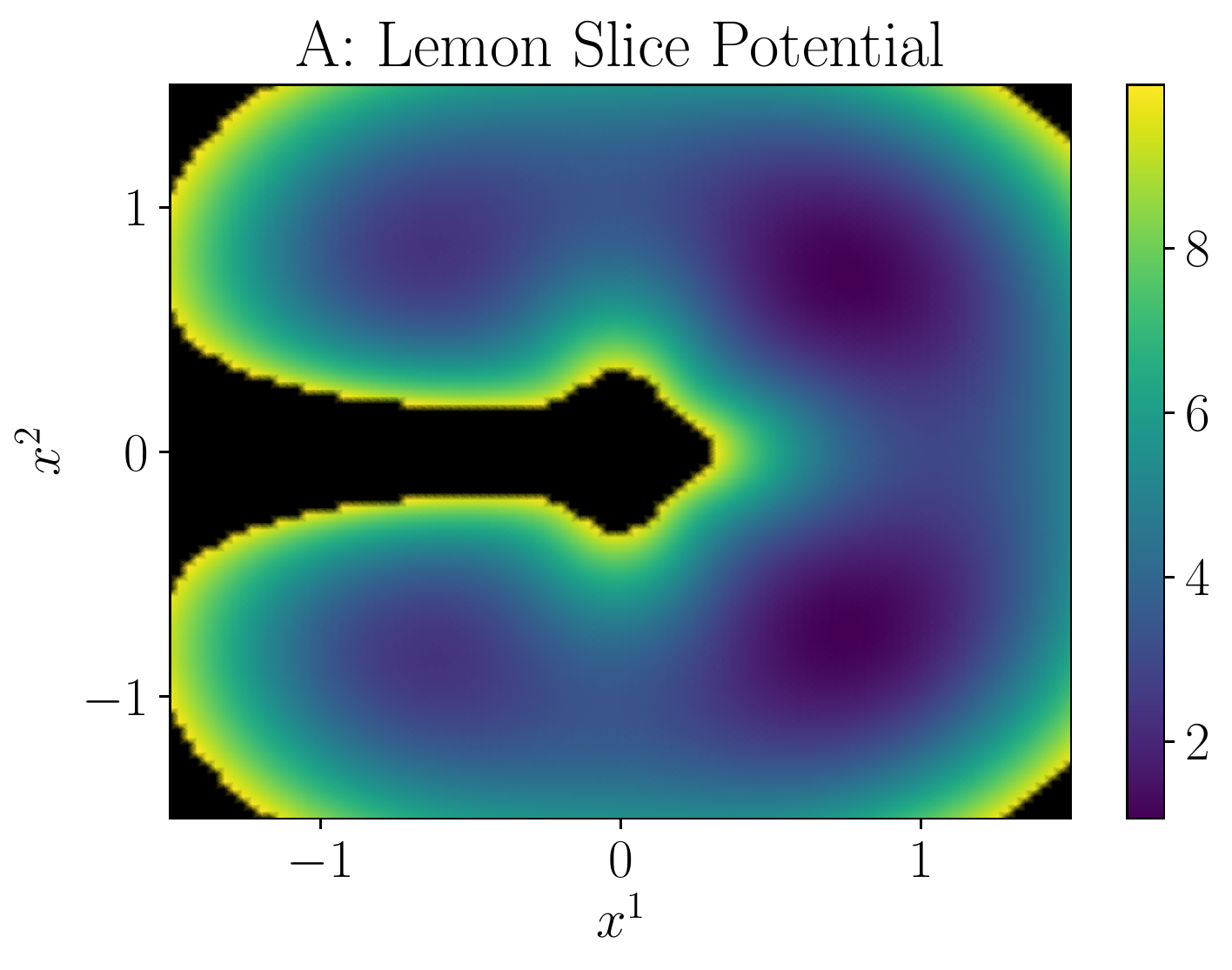}
\includegraphics[width=0.48\textwidth]{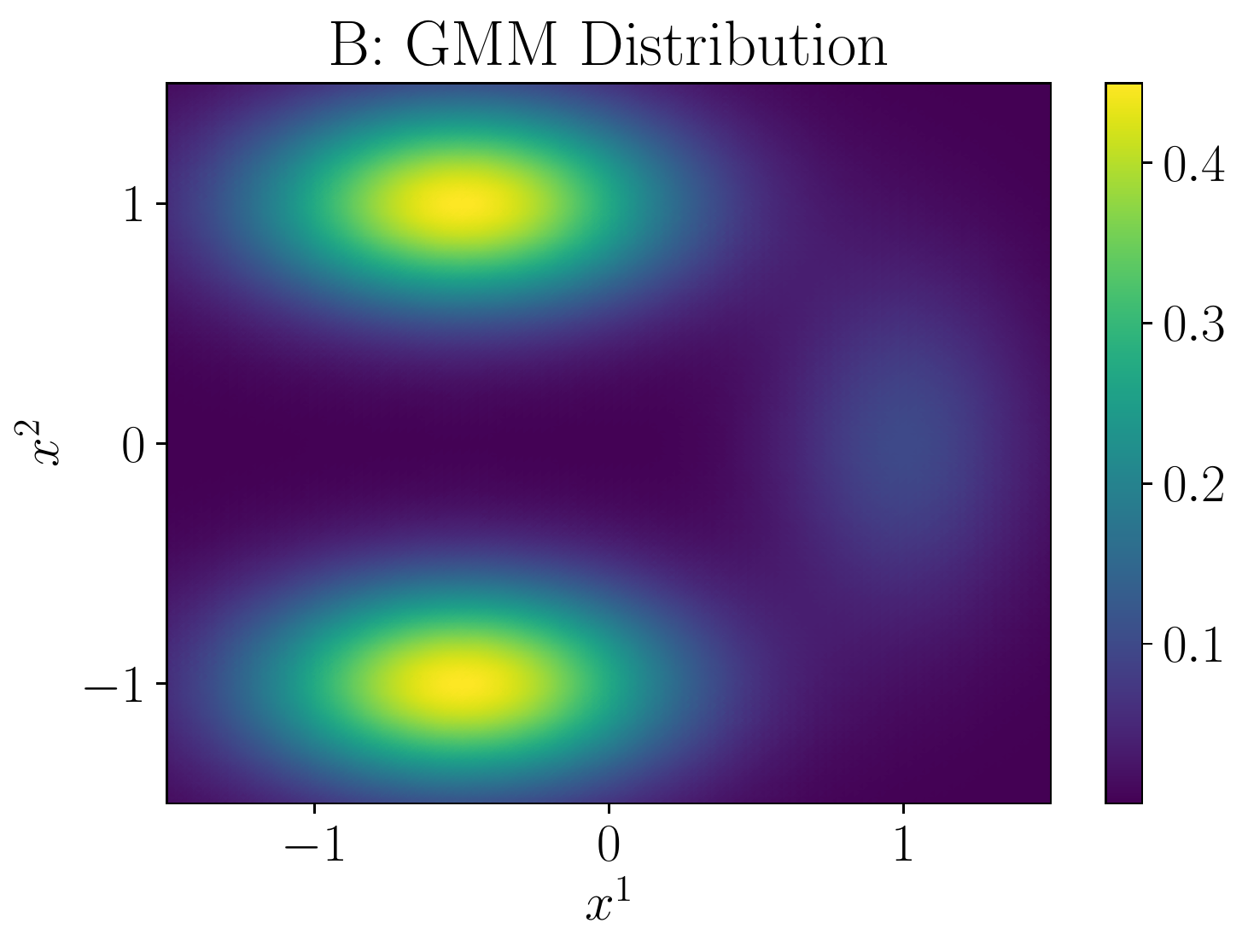}
\caption{A: Contour plot of the two-dimensional Lemon Slice potential Eq.~\eqref{eq:lemon_slice}. B: Contour of the Gaussian mixture distribution of the data used in combination with the re-weighting scheme. \label{fig:lemon_contour}}
\end{figure}

\paragraph{Description} The first system we study is defined by the reversible SDE~\eqref{eq:sde} on $\mathbb{X} = \mathbb{R}^4$, with drift given by the negative gradient $b = -\nabla V$ of a scalar potential $V$ and the diffusion equal to a constant multiple of the identity, $\sigma \equiv \sqrt{2}\mathrm{Id}$. The potential $V$ is a sum of the two-dimensional Lemon Slice potential $V_{\mathrm{LS}}$ along the first two state variables and harmonic potentials $V_h$ along the last two variables, that is
\begin{align}
\label{eq:lemon_slice}
V &= V_{\mathrm{LS}}(x^1, x^2) + V_h(x^3) + V_h(x^4), & \\
V_{\mathrm{LS}}(x^1, x^2) &= \cos(4\varphi) + \frac{1}{\cos(0.5 \varphi)} + 10\ts (r - 1)^2 + \frac{1}{r}, & V^h(x^i) &= 5(x^i)^2.
\end{align}
Above, $r, \varphi$ are two-dimensional polar coordinates derived from $x^1,\,x^2$. A contour plot of the Lemon slice potential, which has also been used in a number of previous publications \cite{Bittracher2018,Nueske2021a}, is shown in Figure~\ref{fig:lemon_contour} A. The dynamics along each of the variables $x^3,\, x^4$ is independent from that along the remaining ones and equilibrates quickly. The dynamics along $x^1,\, x^2$ on the other hand, is metastable with four long-lived states corresponding to the minima of the potential $V_{\mathrm{LS}}$. Consequently, a Markov state model analysis using a box discretization in $(x^1, \, x^2)$-space yields an estimate of these eigenvalues as $\kappa_0 = 0, \, \kappa_1 = -0.42,\, \kappa_2 = -1.19, \, \kappa_3 = -1.54$, which translate into implied timescales $t_1 = 2.38,\, t_2 = 0.84, \, t_3 = 0.65$. This simple example mainly serves to illustrate the importance of choosing TT ranks appropriately, and also to illustrate the use of the re-weighting technique described in Section~\ref{subsec:extensions}.
\paragraph{Methods} For tgEDMD, we generate ten independent long simulations of the SDE~\eqref{eq:sde} by the Euler-Maruyama scheme at integration time step $\Delta_t = 10^{-3}$, each spanning $3\cdot 10^5$ time steps. These simulations are downsampled to $m = 3000$ data points each. Moreover, we generate ten independent data sets of the same size, drawn from a three component Gaussian mixture distribution (GMM) shown in Figure~\ref{fig:lemon_contour} B, in order to test the re-weighting method. The parameters of the GMM were chosen in such a way that samples will correspond to physically relevant states with high probability, while drastically shifting the relative weights of the four metastable states.

\noindent We define elementary subspaces $\mathbb{V}^k, \, 1\leq k \leq 4$, each spanned by Gaussian basis functions depending on $x^k$ alone. For the first two variables, we use seven Gaussians centered equidistantly in $[-1.2, 1.2]$, each with bandwidth $\sqrt{0.4}$, while for the remaining variables we only use five Gaussians centered in $[-1.0, 1.0]$ with bandwidth $\sqrt{0.5}$. For each data set, we calculate the reduced matrix $\hat{M}_{r_p}^{\mathrm{rev}}$ by Algorithm~\ref{alg:tgedmd}. The truncation threshold for the global SVD is set to $\sqrt{m}\epsilon$, where $\epsilon$ varies between $\epsilon = 10^{-8}$ and $\epsilon = 10^{-2}$. We calculate the first four dominant eigenvalues and implied time scales of the reduced matrix, along with their eigenvectors and the evaluation of their associated eigenfunctions at all data sites. Finally, these eigenfunction trajectories are passed to the PCCA method \cite{Deuflhard2005} to generate a clustering of the data into four metastable states. As a comparison, we also apply the standard AMUSE Algorithm~\ref{alg:AMUSE} to the full product basis, using the same truncation parameters. We record the resulting ranks and implied time scales.

\begin{figure}
\centering
\includegraphics[width=0.48\textwidth]{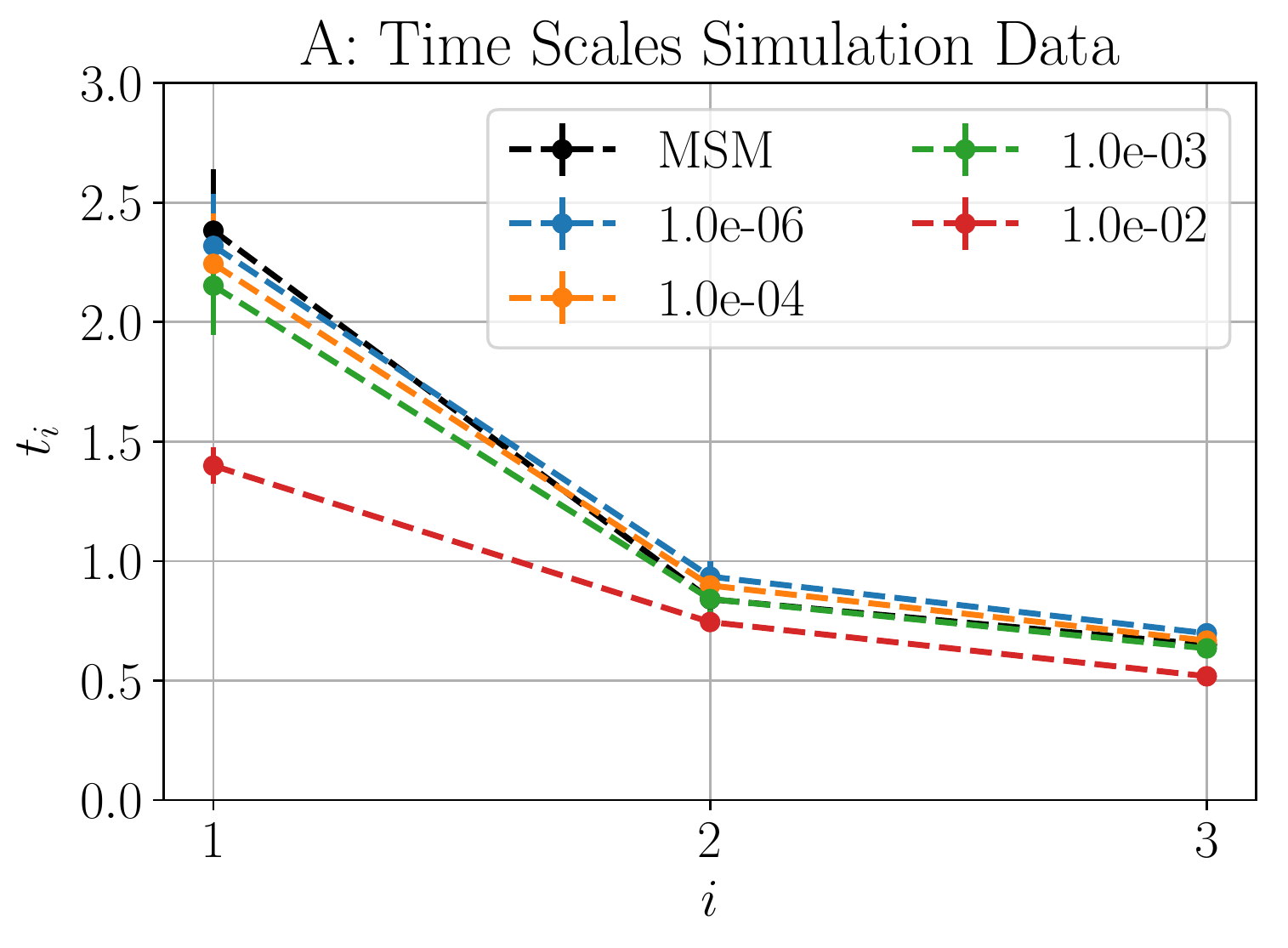}
\includegraphics[width=0.48\textwidth]{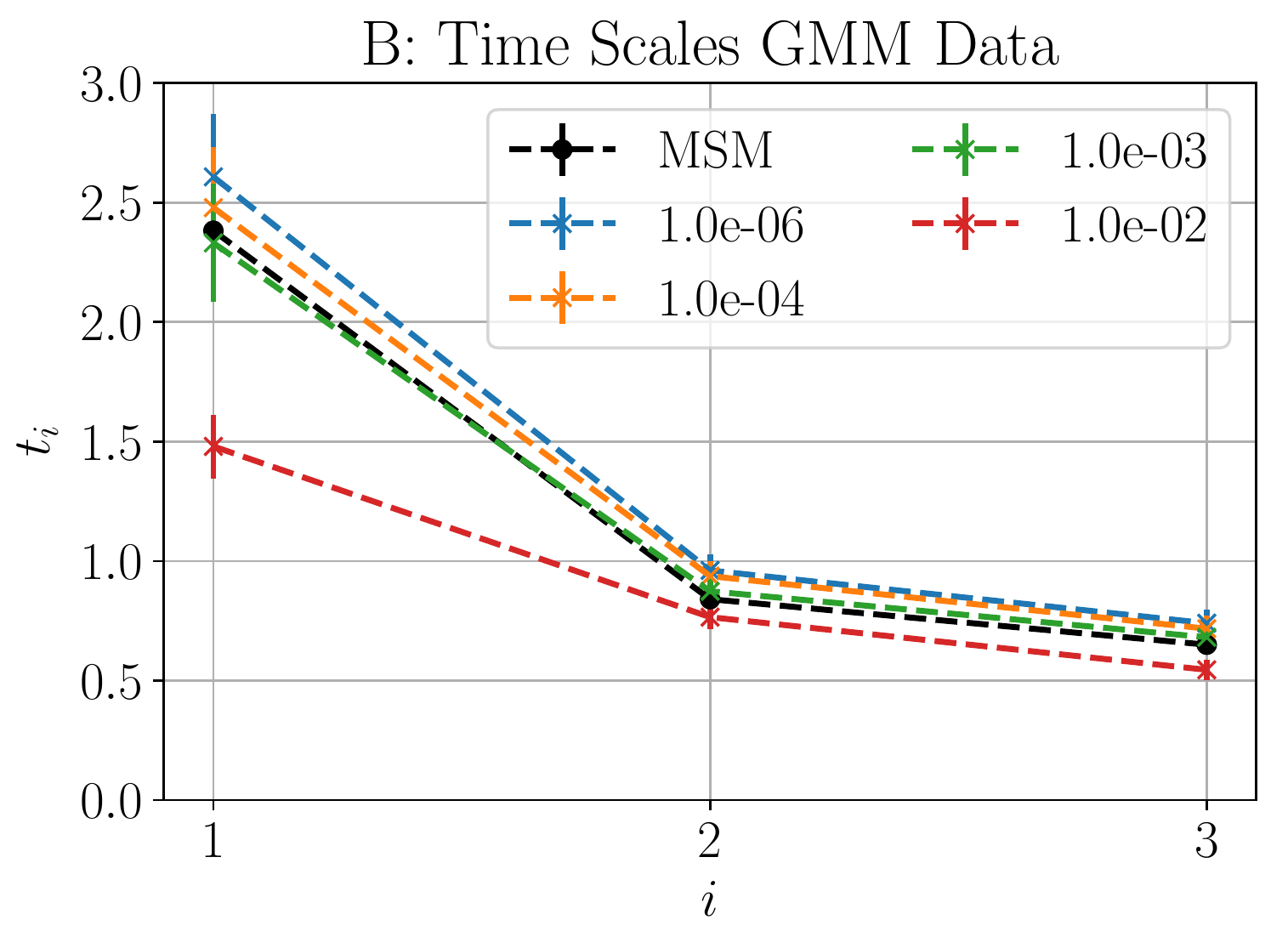}
\includegraphics[width=0.48\textwidth]{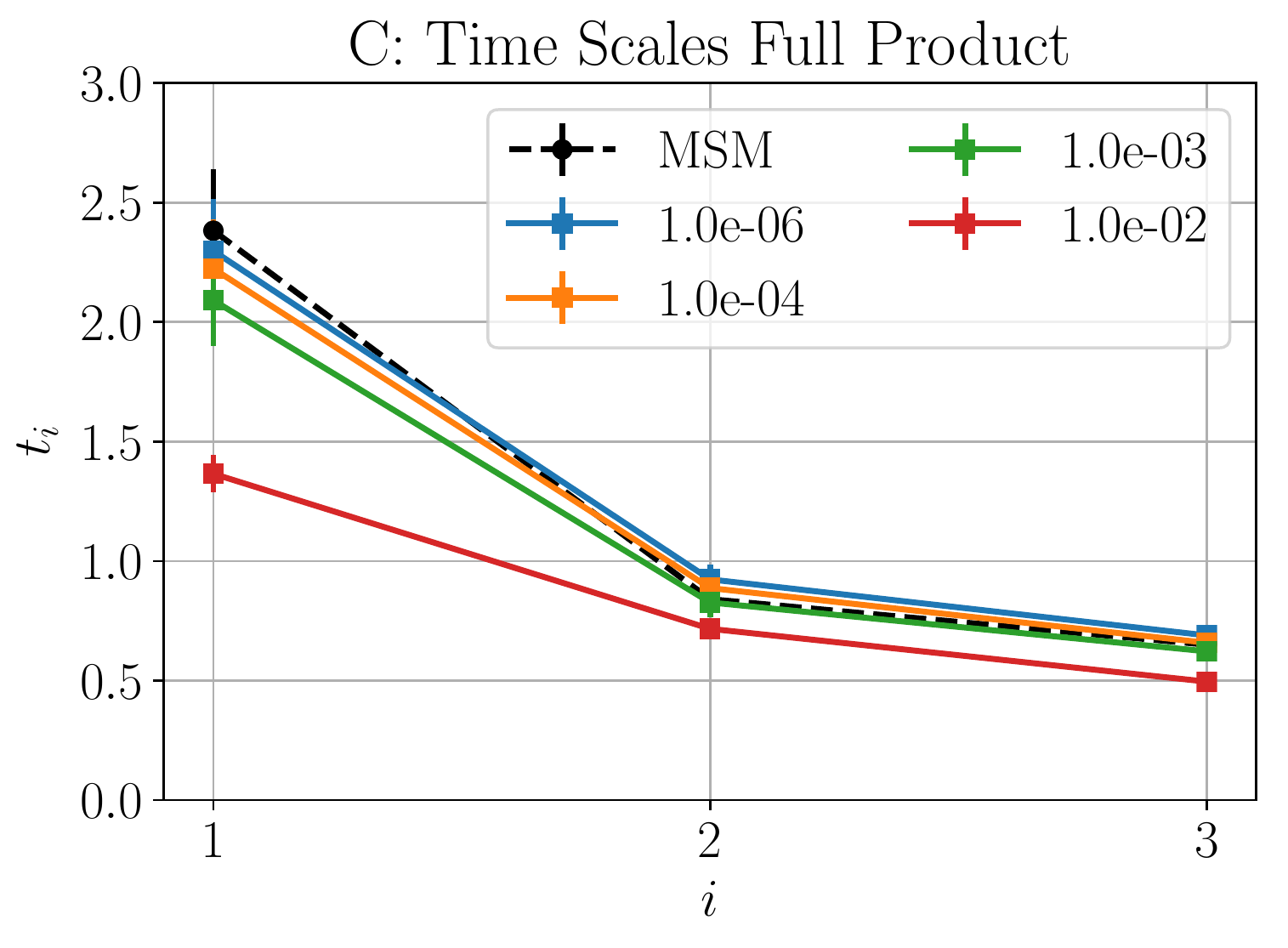}
\includegraphics[width=0.48\textwidth]{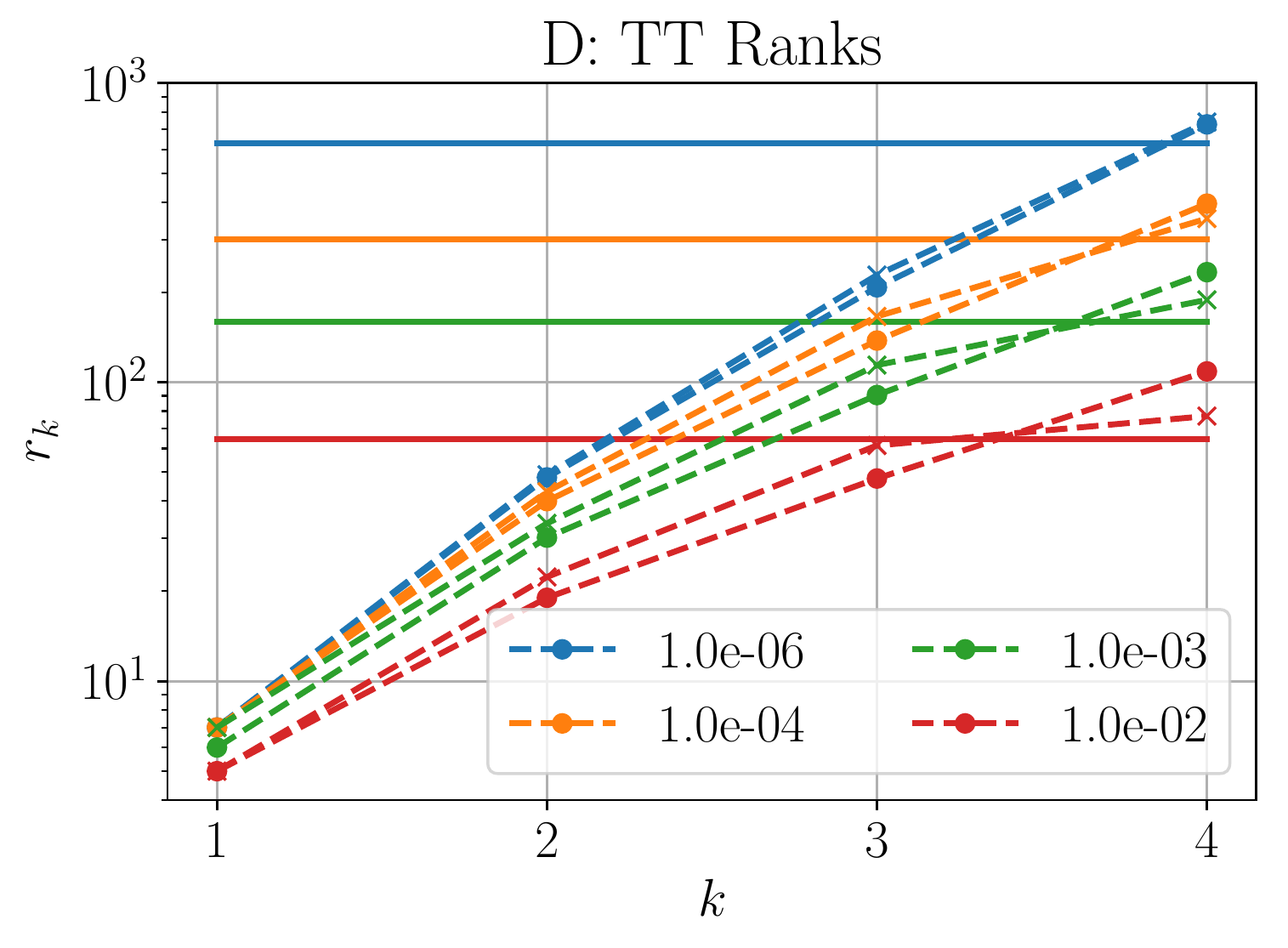}
\includegraphics[width=0.48\textwidth]{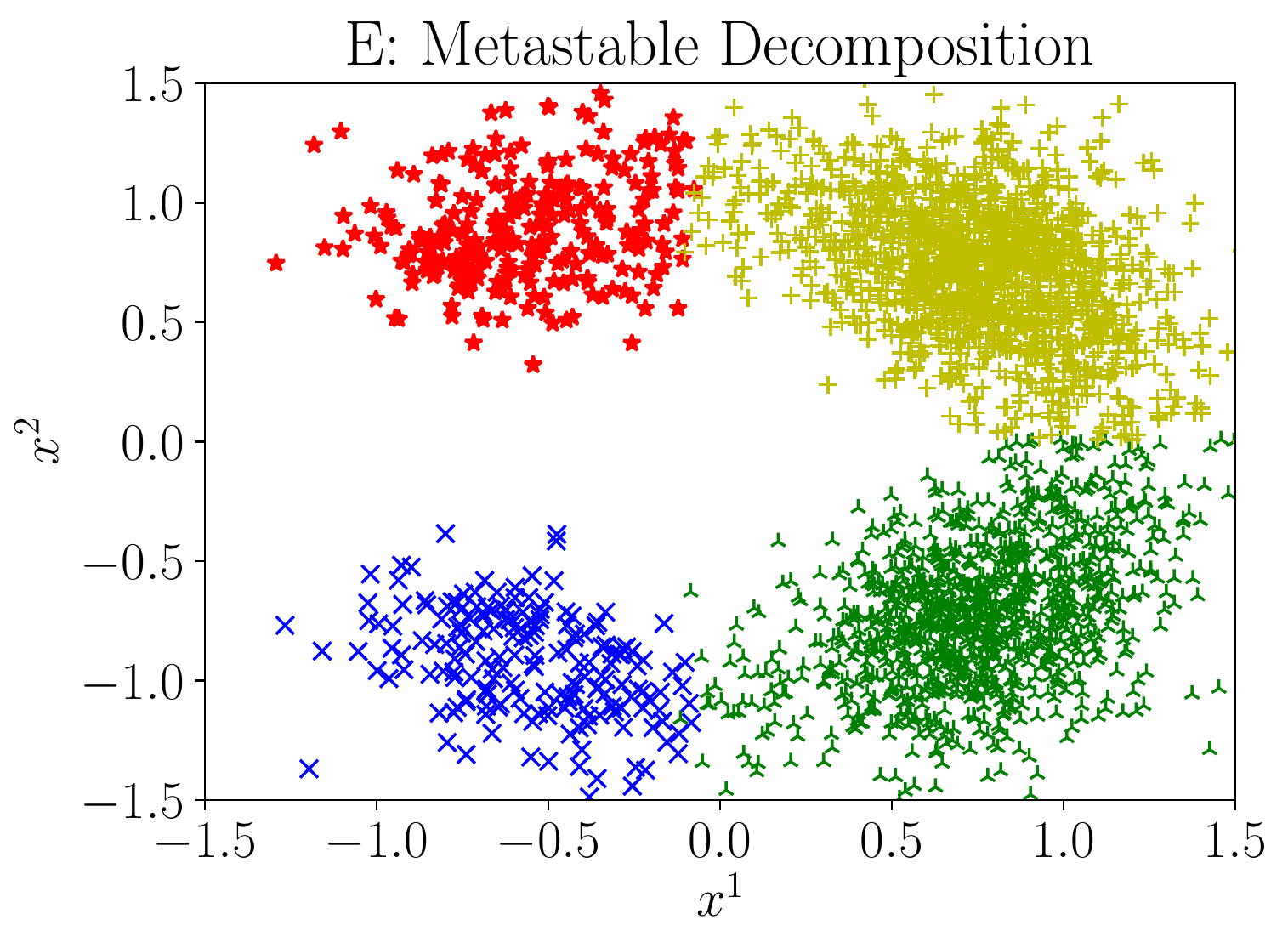}
\includegraphics[width=0.48\textwidth]{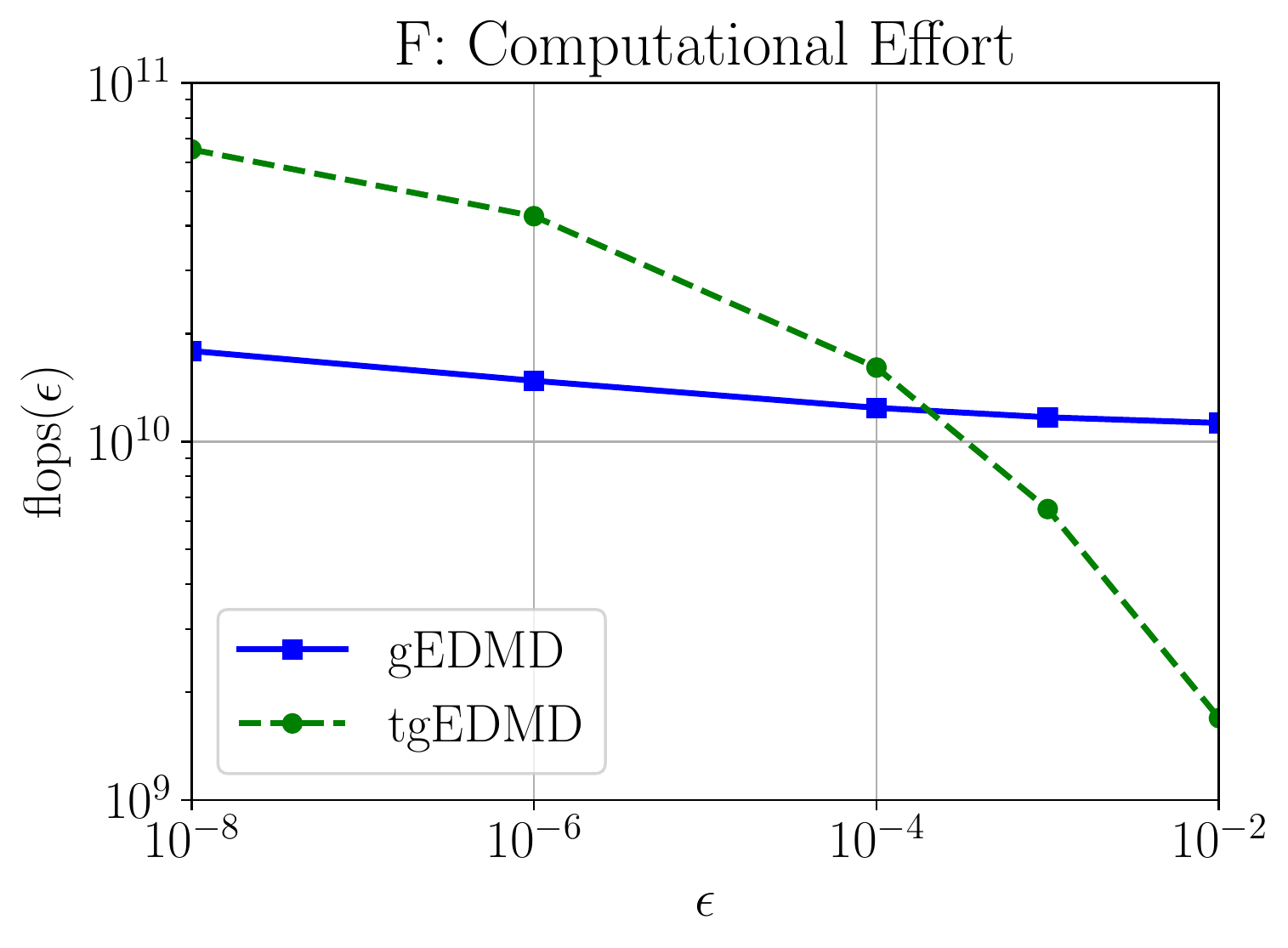}
\caption{Results for four-dimensional diffusion subject to the Lemon Slice potential Eq.~\eqref{eq:lemon_slice} A: Estimates of the leading implied time scales $t_i, \, i=1, 2, 3$, obtained by applying the tgEDMD Algorithm~\ref{alg:tgedmd} with Gaussian basis functions, using different values of the truncation threshold $\epsilon$. MSM-based reference values are shown in black, errorbars were computed over ten independent simulations. B: The same if the data are generated using a three component Gaussian mixture model, combined with the re-weighting method outlined in Section~\ref{subsec:extensions}. C: The same if time scales are calculated using the standard AMUSE Algorithm~\ref{alg:AMUSE} and the SDE data. D: TT ranks $r_k, \, k=1,\ldots, 4$, for different values of the truncation threshold $\epsilon$. Circles indicate results from the SDE trajectories, crosses those obtained using the GMM data. Solid horizontal lines indicate the rank obtained by applying standard AMUSE to the SDE data with truncation threshold indicated by the color. E: Metastable decomposition projected onto $(x^1, x^2)$, obtained by applying PCCA to the eigenfunction trajectories for one of the ten SDE trajectories and $\epsilon = 10^{-3}$. F: Computational cost for tgEDMD and standard gEDMD as estimated by Eqs.~(\ref{eq:comp_cost_tgedmd}-\ref{eq:comp_cost_gedmd}). \label{fig:lemon_slice_results}}
\end{figure}

\paragraph{Results} In Figure~\ref{fig:lemon_slice_results}, we show estimates for the leading three implied time scales $t_i,\, i=1, 2, 3$, obtained from the eigenvalues of the reduced matrix, averaged over ten independent data sets, using both simulation data of the SDE~\eqref{eq:sde} (panel A) and data sampled from the GMM (panel B). Time scales obtained by the standard AMUSE algorithm are shown in panel C as a comparison. We observe that using both the simulation data and the GMM data, we can robustly identify all three time scales within statistical uncertainty by the tgEDMD method, for a range of truncation parameters between $\epsilon = 10^{-6}$ and $\epsilon = 10^{-3}$. If a larger threshold is used, approximation quality starts deteriorating, especially for the slowest time scale. We also verify in panel E that a PCCA analysis of the eigenfunction trajectories computed by tgEDMD yields the correct decomposition of the data set into four metastable states. As shown in panel D of Figure~\ref{fig:lemon_slice_results}, the maximal rank employed by tgEDMD varies between roughly 200 for $\epsilon = 10^{-3}$ and about 800 for $\epsilon = 10^{-6}$. These ranks are approximately equal to the truncation rank achieved by standard AMUSE if the same thresholds are used. In summary, these results confirm that tgEDMD can robustly identify suitable reduced subspaces for spectral analysis of the Koopman generator. Due to the limited size of this example, there is no significant gain in terms of computational effort however. In fact, we see in panel F that estimators Eqs.~(\ref{eq:comp_cost_tgedmd}-\ref{eq:comp_cost_gedmd}) yield about the same computational cost for tgEDMD and standard gEDMD for all truncation ranks. One can observe, however, that the cost for tgEDMD depends more strongly on the truncation rank, whereas the cost of gEDMD is dominated by the basis set size $N = \prod_{k=1}^p n_k$, as expected.

\subsection{Deca-Alanine}
\label{subsec:deca_alanine}

\begin{figure}
	\centering
	\includegraphics[width=0.45\textwidth]{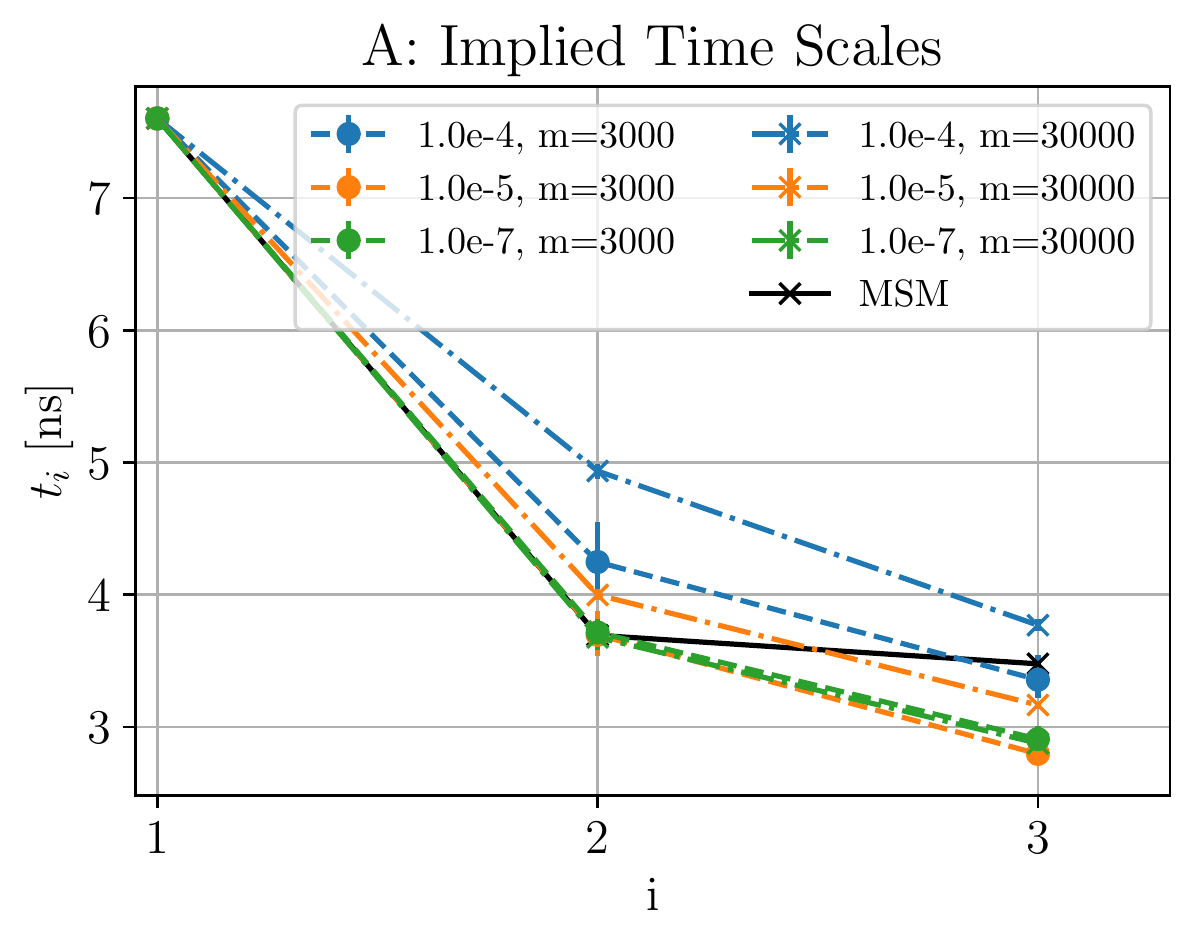}
	\includegraphics[width=0.45\textwidth]{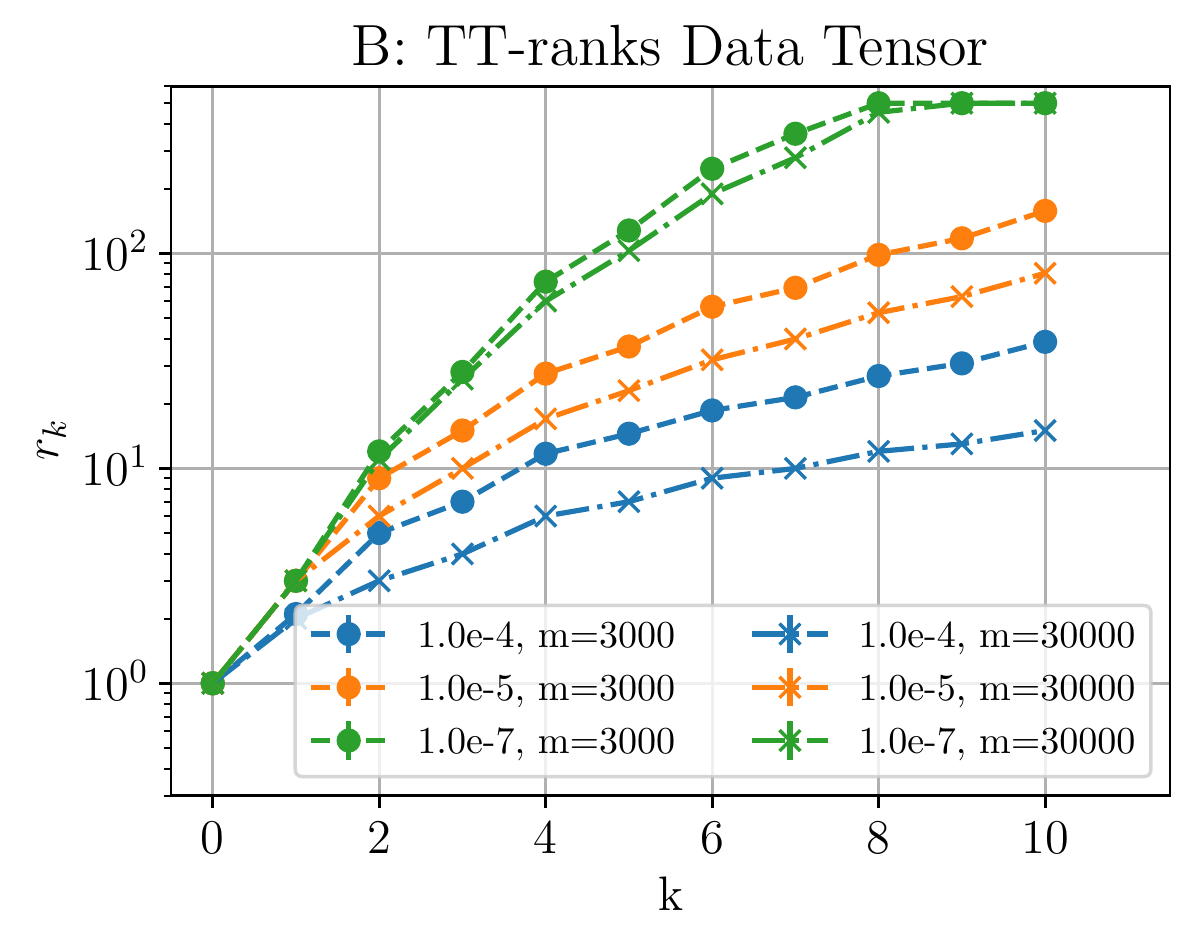}
	\makebox[\textwidth][r]{\includegraphics[width=1.0085\textwidth]{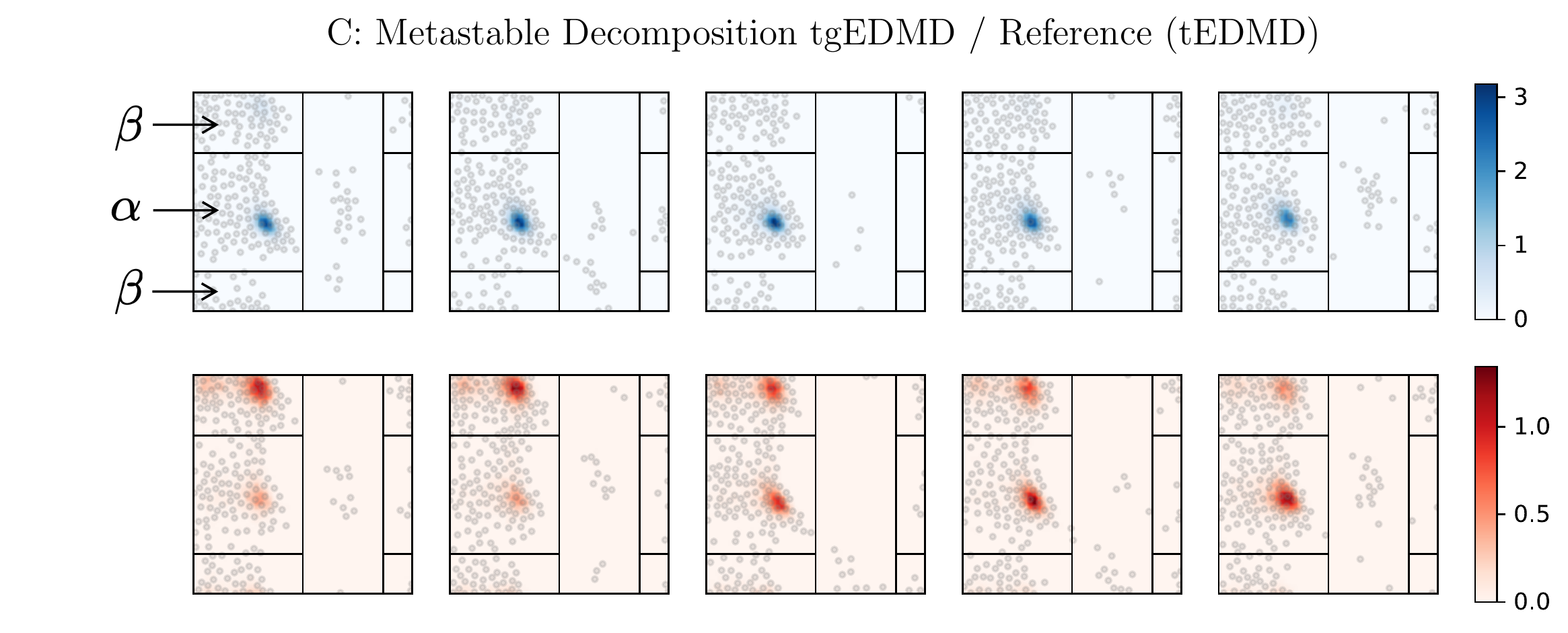}}
	\makebox[\textwidth][r]{\includegraphics[width=0.90\textwidth]{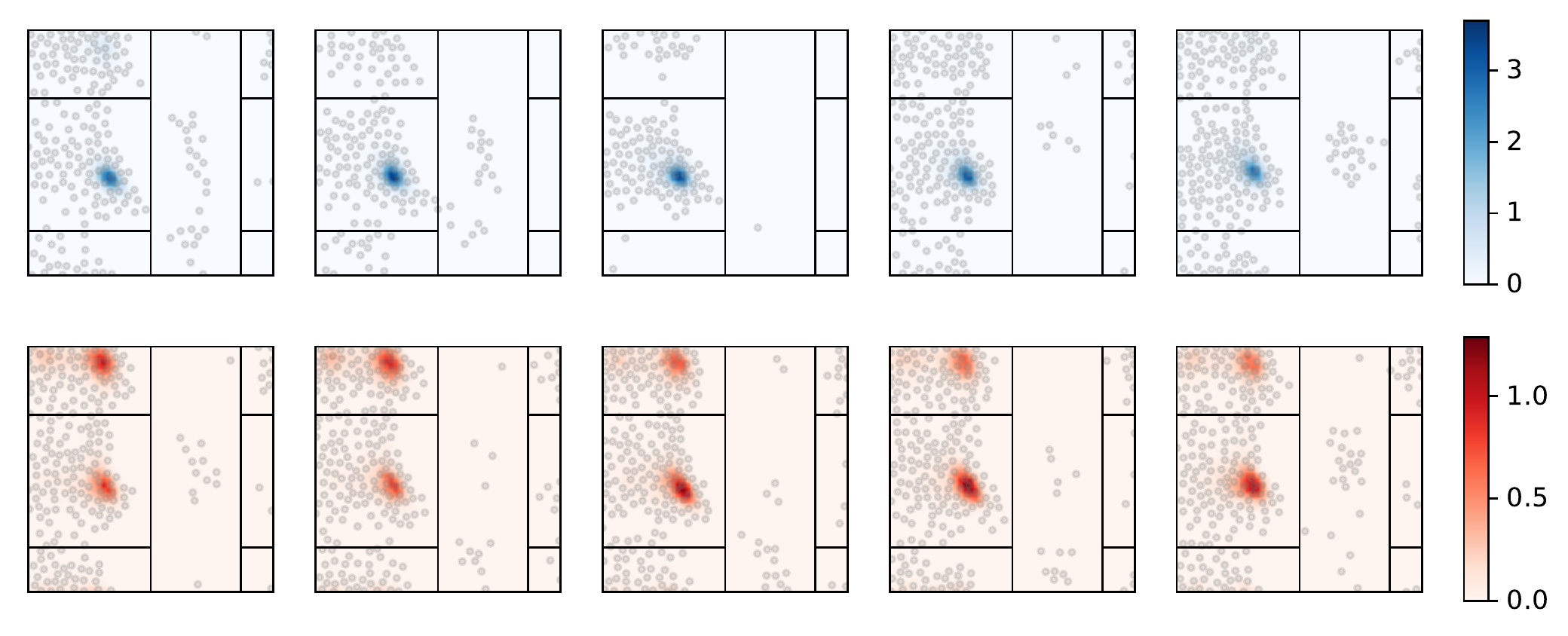}}
	\caption{Results for deca-alanine peptide.
		A: Estimates of the leading three implied time scales, obtained by applying the tgEDMD Algorithm~\ref{alg:tgedmd} with Gaussian basis functions, using different values of the truncation threshold $\epsilon$ and different data sizes. Errorbars were generated by selecting ten different subsets of the total simulation data set, which comprises $3\cdot 10^5$ points. MSM-based reference values are shown in black.
		B: TT-ranks for different values of the truncation threshold $\epsilon$ and different data sizes.
		C: Metastable decomposition of the state space, obtained by applying PCCA to the tgEDMD eigenfunction trajectories ($m=30000$ and $\epsilon = 10^{-7}$). First row: density of data points assigned to the first PCCA state (folded, blue) along the five Ramachandran planes. Second row: the same for the second PCCA state (unfolded, red).
		Third and fourth row: densities of points assigned to the folded and unfolded states by applying PCCA to the eigenfunctions obtained by tensor-based EDMD as in~\cite{Nueske2021}. \label{fig:ala_results}}
\end{figure}

\paragraph{Description}
The second example is a data set resulting from equilibrium molecular dynamics simulations of the ten residue peptide deca-alanine. The simulation setup can be found in \cite{NUESKE2016}, we only note that the data include both positions and momenta of all atoms, and the velocity re-scaling thermostat was used to control temperature.
The (downsampled) data set consists of $m=3 \cdot 10^5$ frames at a time spacing of $10\ \text{ps}$. The system has been analyzed in several previous publications, e.g.  \cite{Nueske2014,Vitalini2015,NUESKE2016,Nueske2021}. It is well-known that the system possesses two major metastable conformations, namely the folded (helical) and unfolded state. These can be distinguished by a collective transition in the space of backbone dihedral angles $\varphi, \psi$. Representing each $(\varphi,\psi)$-pair by its Ramachandran plane, the folded state is (approximately) characterized by $\varphi \leq 0^\circ \lor \varphi \geq 130^\circ$ and $\psi \in [-115^\circ, 80^\circ]$. We will refer to this region in the Ramachandran plane as the $\alpha$-region. In the unfolded state we expect to see a mixture of angle pairs in the $\alpha$-region and in the $\beta$-region, which is defined by $\varphi \leq 0^\circ \lor \varphi \geq 130^\circ$ and $\psi \leq -115^\circ \lor \psi \geq 80^\circ$ (see the upper left plot of Figure~\ref{fig:ala_results} C). The slowest implied time scale, corresponding to the folding process, was determined as $t_1 \approx 7.5\,\mathrm{ns}$ by an MSM analysis in~\cite{Nueske2021}. Two additional time scales were estimated as $t_2 \approx 3.7\,\mathrm{ns},\, t_3 \approx 3.4\,\mathrm{ns}$.

\paragraph{Methods}
As the physically relevant dynamics of deca-alanine can be described by dihedral angle coordinates, we use elementary subspaces of functions defined on ten internal dihedral angles for tgEDMD, the outermost dihedrals being more flexible. The elementary subspaces $\mathbb{V}^k$ are the same as in~\cite{Nueske2021}, they are spanned by periodic Gaussian functions
\begin{align*}
	\psi_{k, u_k}(x) = \exp\Big[ -\frac{1}{2 s_{k, u_k}} \sin^2(0.5(x_k - c_{k, u_k})) \Big].
\end{align*}
Note that $\psi_{k, u_k}(x)$ only depends on the $k$-th coordinate of the state $x$.
If the $k$-th coordinate refers to a $\varphi$ angle, we choose the basis set comprised of the constant function and two periodic Gaussians with $s_{k, u_k} = 0.8, 0.5$ and $ c_{k, u_k} = -2.0,1.0$.
If the $k$-th coordinate refers to a $\psi$ angle, we choose the basis set as the constant function and three periodic Gaussians with $s_{k, u_k} = 0.8, 4.0, 0.8$ and $ c_{k, u_k} = -0.5, 0.0, 2.0$.
Thus, the full tensor space $\mathbb{V}$ consists of $3^5 \cdot  4^5 \approx 2.5 \cdot 10^5$ basis functions.

\noindent As the basis set is defined in terms of reduced coordinates, the effect of this projection needs to be taken into account as described in Section~\ref{subsec:extensions}. The full state generator is then, in principle, given as the generator of the molecular dynamics engine on full position and momentum space. For several variants of thermostatted MD, including Langevin dynamics and velocity re-scaling, it is well-known that the dynamics on the atomic position space can be approximated by a reversible SDE with constant diffusion on long time scales (see~\cite{Nueske2021a} for a detailed discussion). However, this approximation induces a re-scaling of time by a typically unknown factor. Still, we make use of this approximation here by employing the reversible tgEDMD algorithm. Besides the basis functions and their derivatives, the only additional quantity we require is the Jacobian of all dihedral angles with respect to the atomic positions, which can be calculated analytically.

\noindent We then compute the reduced matrix $\hat{M}_{r_p}^{\mathrm{rev}}$ by Algorithm~\ref{alg:tgedmd} using downsampled trajectories of either $m = 3000$ or $m = 30000$ snapshots.
In this example, the dominant singular values of the TT-cores are rather large and differ significantly in size between the cores. Therefore, we employ a relative SVD truncation threshold in contrast to the absolute cut-off used before, i.e. the truncation threshold is set to $\sqrt{m}\epsilon s_1$, where $s_1$ denotes the largest singular value, and $\epsilon$ varies between $\epsilon=10^{-7}$ and $\epsilon=10^{-4}$. In addition, we cap the maximally allowed SVD rank at $500$.
As before, we calculate the dominant eigenvalues and implied time scales of the reduced matrix. In order to account for the change of time units, we calculate the ratio of the reference slowest time scale $t_1$ and our estimate and re-scale all time scales by that ratio. Of course, this enforces all tgEDMD estimates for $t_1$ to match the reference exactly. However, by also comparing the next few estimated time scales to their reference values, we can verify whether or not our estimates really match the reference up to a re-scaling of time. Finally, we also compute a clustering of the data into metastable states using the PCCA method \cite{Deuflhard2005}.

\begin{figure}
	\centering
	\includegraphics[width=0.7\textwidth]{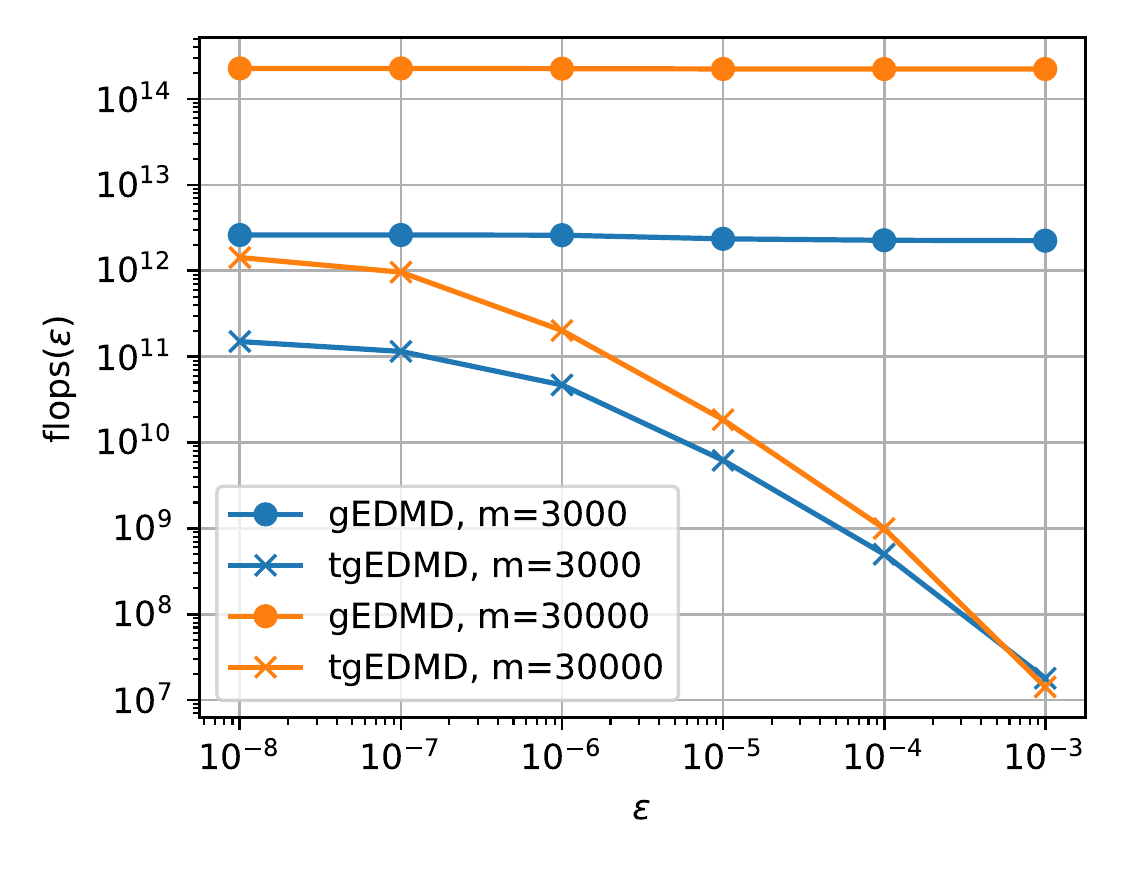}
	\caption{Computational cost of tgEDMD compared to standard AMUSE (gEDMD), as estimated by Eqs.~(\ref{eq:comp_cost_tgedmd}-\ref{eq:comp_cost_gedmd}), for the deca-alanine example.}
	\label{fig:ala_complexity}
\end{figure}

\paragraph{Results}
In panel A of Figure~\ref{fig:ala_results}, we show estimates for the leading three implied time scales $t_i,\, i=1, 2, 3$, obtained from the eigenvalues of the reduced matrix. We find that these time scales can be reliably approximated across the range of parameters and data sizes we tested. A ten to fifteen per cent relative error margin is only exceeded for the largest SVD threshold $\epsilon = 10^{-4}$.\\
The ranks of the associated TT representations of the data tensor are shown in panel B of Figure~\ref{fig:ala_results}. The maximal TT rank of $500$ is only attained for the smallest threshold $\epsilon = 10^{-7}$. The displayed ranks are similar to the results presented in~\cite{Nueske2021} and we also confirm their finding that TT ranks on the order of $50$ to $100$ are sufficient to satisfactorily approximate the slow dynamics of the peptide.  \\
In panel C, we verify that the folded and unfolded states were correctly identified by tgEDMD. To this end, we separately compute histogram densities of all data points that were assigned to each metastable state, projected onto each of the five Ramachandran planes. The results obtained for parameters $m=30000$, $\epsilon=10^{-7}$, are shown in the first two rows. These results are consistently reproduced for all thresholds and data sizes we tested. As a reference, we show the same densities based on the clustering obtained from the eigenfunctions computed by tensor-based EDMD as in~\cite{Nueske2021} (last two rows of panel C).
We observe that both the folded and unfolded states were identified correctly, up to a slight discrepancy of the relative weights of the alpha- and beta-regions of the first two planes for the unfolded state.\\
In this ten-dimensional example, the computational benefit of employing the presented tgEDMD method compared to the standard AMUSE algorithm (gEDMD) is already substantial, as Fig.~\ref{fig:ala_complexity} illustrates. For the parameters we tested, the computational effort of tgEDMD is orders of magnitude smaller than for gEDMD. For instance, we see that by setting $\epsilon = 10^{-5}$ and $m=30000$, accurate results can be achieved (cf. Fig.~\ref{fig:ala_results}) at computational cost reduced by more than three orders of magnitude. Moreover, we confirm once again that the computational effort for tgEDMD depends primarily on the TT ranks, thus offering an efficient means of controlling the computational cost. On the other hand, the cost of gEDMD is dominated by the basis set size $N$ and essentially independent of the truncation rank.

\section{Conclusions}
\label{sec:conclusions}
We have presented tgEDMD, a data-driven method to approximate the Kolmogorov backward operator, or generator, of a dynamical system driven by a stochastic differential equation using the tensor train format. The centerpiece of the method is a TT representation of the generator data tensor. We have derived this representation for reversible and non-reversible SDEs. We have shown that tgEDMD consistently approximates the generator on a multi-linear subspace of the full tensor space and also analyzed the computational complexity. We have shown how importance sampling methods can be incorporated and how the method can be used to compute an effective generator on a reduced state space. Finally, we have presented successful applications of the method to a low-dimensional test system and to a benchmarking molecular dynamics data set.

\section*{Acknowledgments}
The authors thank the Paderborn Center for Parallel Computing (PC2) for computational resources. M. L. has been supported by Deutsche Forschungsgemeinschaft (DFG) under Germany’s Excellence Strategy via the Berlin Mathematics Research Center MATH+ (EXC-2046/1 project ID:  390685689).

{\small{}\bibliographystyle{unsrturl}
\bibliography{Library}
}{\small\par}

\end{document}